\documentclass{jssc}

\usepackage[toc,page,title,titletoc,header]{appendix}
\usepackage{dsfont}
\usepackage{color}
\makeatletter
\def\@oddfoot{\hfill}
\newcount\shujiaocount
\def\setshujiao{%
  \shujiaocount=\count0
  \def\@oddfoot{%
      \ifodd\count0
      \else
      \fi
  }%
}
\makeatother

\def\ba{\begin{array}}
\def\ea{\end{array}}
\def\ban{\begin{eqnarray*}}
\def\ean{\end{eqnarray*}}
\def\be{\begin{equation}}
\def\ee{\end{equation}}
\def\bna{\begin{eqnarray}}
\def\ena{\end{eqnarray}}

\def\nn{\nonumber}

\def\dref#1{(\ref{#1})}

\def\textsubscript#1%
{$_{\text{#1}}$}
\newcommand*\supercite[1]{\textsuperscript{\cite{#1}}}
\def\cdd{\mbox{\boldmath$\cdot$}~}

\def\ee{{\rm e}}
\newcommand{\rulex}{\hfill\rule{1mm}{3mm}}
\def\ay{\arraycolsep=1.5pt}
\def\d{\displaystyle}
\def\dfrac{\displaystyle\frac}
\input{amssym.def}
\include{graphix}

\usepackage[font=small,labelsep=quad,labelfont=bf,width=\textwidth,format=hang,indention=0mm]{caption}

\makeatletter
\def\@oddfoot{\hfill}
\newcount\shumeicount
\def\setshumei#1#2#3{%
  \shumeicount=\count0
  \def\@oddhead{%
    \raise-5pt\hbox to0pt{\vrule width\hsize height 0pt depth 0.4pt\hss}\relax
    \ifnum \shumeicount=\count0
      \raise-7pt\hbox to0pt{\vrule width\hsize height 0pt depth 0.4pt\hss}\relax
      #1
    \else
      \ifodd\count0
        #2
      \else
        #3
       \fi
     \fi
  }%
}
\makeatother
\makeatletter
\def\@oddfoot{\hfill}
\newcount\shujiaocount
\def\setshujiao{%
  \shujiaocount=\count0
  \def\@oddfoot{%
      \ifodd\count0
      \else
      \fi
  }%
}
\makeatother
\def\biaoti#1#2#3#4{{
  \vspace*{0.3cm}
  \begin{flushleft} \Large\bf #1\end{flushleft}
  \vspace*{-0.2cm}
      \begin{flushleft}
      \bf #2
      \end{flushleft}
      \footnotetext{\hspace{-6mm} #3\\ #4}}}

\def\dshm#1#2#3#4
{\setshumei{J Syst Sci Complex (#1) #2:
{\thepage--\pageref{LastPage}}\hfill}
            {\hfill {\small #3}\hfill\hbox to0pt{\hss\thepage}}
            {\hbox to0pt{\thepage\hss}\hfill {\small #4}\hfill
            }
            \setshujiao}
\def\drd#1#2
{{\vskip 1cm\small \begin{flushleft}
 #1 \\
 #2\\
\copyright The Editorial Office of JSSC \&  Springer-Verlag GmbH
Germany 2019
\end{flushleft}}}

\def\bar{\overline}
\def\epsilon{\varepsilon}
\def\proof{\vspace{1mm}\indent {\it Proof}\quad}

\begin{document}

\biaoti{Flocking with General Local Interaction and Large Population$^*$}
{\uppercase{CHEN} Ge \cdd \uppercase{LIU}
Zhixin}
{\uppercase{CHEN} Ge \cdd \uppercase{LIU} Zhixin\\
{\it Key Laboratory of Systems and Control $\&$ National Center for
Mathematics and Interdisciplinary Sciences, Academy of Mathematics and Systems Science,
 Chinese Academy of Sciences, Beijing $100190$, China.}\\
    Email: chenge@amss.ac.cn;  Lzx@amss.ac.cn.} 
{$^*$The research was supported by the National Natural Science
Foundation of China under Grant No. 11688101, 91634203, 91427304, and 61673373, and
 the National Key Basic Research Program of China (973
Program) under Grant No. 2016YFB0800404.
Part of this paper was presented in \emph{the $10$th World Congress on Intelligent Control and Automation},
 pp. 3515--3519, July 6--8, 2012, Beijing, China.
}




\Abstract{This paper studies a flocking model in which the
interaction between agents is described by a general local nonlinear
function depending on the distance between agents. The existing
analysis provided sufficient conditions for flocking under an
assumption imposed on the system's closed-loop states; however this
assumption is hard to verify. To avoid this kind of assumption the
authors introduce some new methods including large deviations theory
and estimation of spectral radius of random geometric graphs. For
uniformly and independently distributed initial states, the authors
establish sufficient conditions and necessary conditions for
flocking with large population. The results reveal that under some
conditions, the critical interaction radius for flocking is almost
the same
 as the critical radius for connectivity of the initial neighbor graph.}      




\section{Introduction}\label{intro}

The flocking, which means the collective coherent motion of a large
number of self-propelled organisms, is an amazing phenomenon in
nature. It attracts the researchers from diverse fields, including
biology, physics, computer science, mathematics and control theory,
see \cite{Toner,Reynolds,Vicsek1,Buhl2006,Chazelle2014,flocking2006}
among many others. In order to investigate the flocking phenomena
exhibited in biological systems, the well-known Vicsek model was
proposed in \cite{Vicsek1}, in which the heading is updated
according to the headings of the corresponding neighbors. Simulation
results reveal that the headings of all agents will be almost the
same (called consensus or flocking) for large population size. The
theoretical study for the flocking behavior of the Vicsek model can
be found in \cite{Jad1,sakvin,jiang,Chen2017}.
Inspired by the Vicsek model, Cucker and Smale also developed a
flocking model in which each agent interacts with all other agents
and the interaction weights decay according to the agents' distance,
and proposed some sufficient conditions for flocking which depend on
system parameters only\supercite{Smale2007}. Following their work,
the theoretic results with different scenarios, for example, the
noisy environment, and the hierarchical structure, are given for the
Cucker-Smale model and its variants (see
\cite{Cucker2008,Cucker2011,Peszek2014,Carrillo2010,Ahn2010,Shen2007,Park2010,Ha2010,Ha2017}).

The assumption of global interaction between agents in the
Cucker-Smale model may not be applicable for some practical systems
especially for systems with large population size. For example, in
\cite{Buhl2006}, Buhl, et al. observed that the individual adjusts
its direction to align with its neighbors within an interaction
range; In \cite{Rosenthal2015}, Rosenthal,  et al. found the golden
shiner (a kind of fish) uses simple, robust measures to assess
behavioral changes of the others in neighboring region, and the
interaction weight depends on the metric distance and ranked angular
area; In \cite{Rieu2000}, Rieu, et al. showed that the spatial
correlation of the velocities of  {Hydra} cells decreases to zero
rapidly with regarding the distance between cells, and remains zero
when the distance is large. For some systems, it may be even
impossible to obtain the explicit expression of interaction between
agents (see \cite{Rieu1998}). To be more practical,  Martin, {et
al.} proposed a multi-agent model in which the interaction between
agents is represented by a general non-negative and non-increasing
function\supercite{Martin2014}. Also, they showed that the system
will achieve flocking behavior if the maximal perturbation allowed
on the relative positions of agents is not bigger than a constant,
and the maximum difference among the initial velocities of all
agents is small enough\supercite{Martin2014}. However, how to
guarantee the condition of the relative position perturbation is
unsolved.

This paper aims at providing some flocking conditions depending on
initial state and system  parameters only for a multi-agent model
with general local interaction.
This problem is challenging because the positions and velocities of
all agents are coupled.  Also, the local interaction between agents
are described by a general state-dependent function. These
characteristics make that many traditional tools like Lyapunov
method cannot be used.
Similar to \cite{tang2007,liu2009b,Chen2014}, we investigate the
flocking behavior under random  initial states. The large deviations
theory is introduced to estimate the maximum degree of the dynamic
weighted graphs, and the extension of the random geometric graph
theory is employed to estimate the spectral radius of the initial
neighbor graphs. We establish sufficient conditions and  necessary
conditions for the flocking behavior of our multi-agent model. In
particular, we reveal that under some conditions,  the smallest
possible interaction radius  for flocking is almost the same as the
critical connectivity radius of the initial neighbor graphs.

The rest of this paper is organized as follows. In Section
\ref{prepa}, we provide the problem formluation. The main results
are given in Section \ref{mainresult}.  In Section
\ref{Sec_estimation}, we present the approach to estimating  and
calculating a key value which characterizes the property of
nonlinear functions. The proofs of our main results are put in
Section \ref{Sec_proof_mrs}. Concluding remarks are made in
Section{}\,\ref{conclusion}.

\section{Problem Statement}\label{prepa}

\subsection{A Nonlinear Flocking Model}
This paper considers a nonlinear discrete-time multi-agent system
composed of $n$ spatially distributed agents, each moving in a
$d$-dimensional $(d\geq 2)$ Euclidean space $\mathbb{R}^d$. Let
$\mathcal{V}\triangleq\{1, 2,\cdots, n\}$  to be the set of all
agents. Set $X_i(t)\in \mathbb{R}^d$ and
 $V_i(t)\in \mathbb{R}^d$ to be the  position and
velocity of agent $i$ at time $t$ respectively. Following \cite{Martin2014} but with some modifications,
for any time $t\geq 0$ and agent $i\in\mathcal{V}$, $X_i(t)$ and $V_i(t)$ are
updated according to the following equation, \bna \label{m1}
\left\{%
\begin{array}{ll}
X_i(t+1)=X_i(t)+V_i(t),\\
\d V_i(t+1)=V_i(t)+\sum_{j=1}^n
f_n\left(\|x_i(t)-x_j(t)\|\right)\left(V_j(t)-V_i(t)\right),
\end{array}%
\right. \ena where the nonlinear function $f_n(\cdot)$ denotes the
local interaction weight, and $\|\cdot\|$ represents the Euclidean
norm. From an intuitive point of view, the larger the distance
between agents, the weaker the interaction weight $f_n(\cdot)$
should be. Without loss of generality, we assume that $f_n(\cdot)$
is a  {non-increasing}  {integrable}  {function} satisfying
\bna\label{m0}
f_n(x)\left\{%
\begin{array}{ll}
>0, ~~~~~~~~0\leq x< r_n,\\
=0, ~~~~~~~~x\geq r_n,
\end{array}%
\right. \ena where $r_n$ is called the  {interaction radius}. In
this paper, the commonly used circular or spherical neighborhood is
adopted. The pair of two agents is called neighbors if and only if
their distance is less than a pre-defined radius $r_n$.

%

Let $x_{ij}(t)$ and $v_{ij}(t)$ denote the $j$th element of $X_i(t)$
and $V_i(t)$, respectively. Denote $X(t):=(x_{ij}(t))_{n\times d}$
and $V(t):=(v_{ij}(t))_{n\times d}$. The system \dref{m1} can be
rewritten into the following matrix form:
\begin{eqnarray}\label{m2}
\left\{%
\begin{array}{ll}
X(t+1)=X(t)+V(t),\\
V(t+1)=P(t)V(t),
\end{array}%
\right.
\end{eqnarray}
where $P(t)=(p_{ij}(t))_{n\times n}$ is the weighted average matrix defined by
\begin{eqnarray*}
  \label{m3}
 p_{ij}(t)=\left\{
             \begin{array}{ll}
               f_n\left(\|x_i(t)-x_j(t)\|\right), & \hbox{if}\ j\neq i, \\
              \d 1-\sum_{j=1, j\neq i}^nf_n\left(\|x_i(t)-x_j(t)\|\right), &
               \hbox{otherwise}.
             \end{array}
           \right.
 \end{eqnarray*}
For any $i\in \mathcal{V}$ and $t\ (=0, 1,   \cdots )$, we have
$\sum_{j=1}^n p_{ij}(t)=1$.

The objective of this paper is to investigate the flocking
behavior of the system \dref{m2} with large population. Follow
\cite{Martin2014} we say that the system (\ref{m2}) achieves a
 {flocking behavior} if  the velocities of all agents reach
agreement, that is, \bna\label{m4}
\lim_{t\rightarrow\infty}\max_{1\leq i,j\leq
n}\|V_i(t)-V_j(t)\|=0.\ena

\subsection{Random Geometric Graphs}

Following our previous work \cite{Chen2014}, we assume that the $n$
agents are independently and uniformly distributed in the unit cube
$[0,1]^d$. Denote $\mathcal{X}_n=\{X_1(0),$ $X_2(0), \cdots,
X_n(0)\}$ as the set of initial positions of the $n$ agents. We
introduce a  {random geometric graph} $G(\mathcal{X}_n;r_n)$ to
describe the neighbor relations between agents at the initial time,
with vertex set $\mathcal{V}$ and with undirected edges connecting
the pairs $\{X_i(0), X_j(0)\}$ satisfying $\|X_i(0)- X_j(0)\|\leq
r_n$; See \cite{penrose2003} for more properties of random geometric
graphs. It is worth pointing out   that the positions are not
independent and the properties of random geometric graphs cannot be
used any more when the agents move around.

In order to state our results clearly, we assume that the interaction radius has
the following expression,
\begin{eqnarray}\label{rad_con}\lim_{n\rightarrow\infty}\left(n r_n^d/\log n\right)= \alpha \in
(0,\infty].
\end{eqnarray} Set $R_c=R_c(n):=\sqrt[d]{\frac{2^{d-1}\log n}{d \pi_d n}}$,
where $\pi_d$ denotes the volume of the unit ball in $\mathbb{R}^d$.
It is proved that for $d=2$, $R_c$ is the  {critical connectivity
radius} of $G(\mathcal{X}_n;r_n)$ in a probability sense
(see{}\,\cite{kumar2000}). The following lemma shows that a similar
result holds for $d\geq 3$.

\begin{lemma}\label{rgg}
{\it The random geometric graph $G(\mathcal{X}_n;r_n)$ is connected
with high probability if $\alpha>\frac{2^{d-1}}{d\pi_d}$, and is not
connected with high probability if $\alpha<\frac{2^{d-1}}{d\pi_d}$,
where $\alpha$ is defined by $(\ref{rad_con})$.}
\end{lemma}

The proof of Lemma \ref{rgg} is presented in Appendix \ref{App_A}.

We call $r_n$ the  {super-critical connectivity radius} of
$G(\mathcal{X}_n;r_n)$ if $\alpha>2^{d-1}/(d\pi_d)$, and we call
$r_n$ the  {sub-critical connectivity radius} of
$G(\mathcal{X}_n;r_n)$ if $\alpha<2^{d-1}/(d\pi_d)$.

\subsection{Large Deviations Techniques}

The large deviations techniques are applied to deal with the influence of
the nonlinear interaction weights $f_n(\cdot)$. We
first introduce some notations.

For a given constant $\delta\geq 0$, define
\begin{eqnarray*}
f_{n,\delta}(x):=\left\{%
\begin{array}{ll}
f_n(0), &\mbox{if $x\leq \delta r_n$},\\
f_n(x-\delta r_n), &\mbox{else}.
\end{array}%
\right.
\end{eqnarray*}
By the definition of $f_{n,\delta}$, we have $f_n=f_{n,0}$. Let
$x_0=(\frac{1}{2},\frac{1}{2},\cdots,\frac{1}{2})$ be the center
point of $[0,1]^d$, and set
$$\xi_{n,\delta}:=f_{n,\delta}\left(\big\|X-x_0\big\|\right),$$
where $X$ is a random variable uniformly distributed in $[0,1]^d$.
For $x\in \mathbb{R}$, define
\begin{eqnarray}\label{com1}
I_{n,\delta}(x):=\sup_{\theta>0}\left\{\theta x-(n-1) \log
\left(E\left[{\rm e}^{\theta\xi_{n,\delta} } \right]\right)
\right\}.
\end{eqnarray}
The function \dref{com1} is called a  {rate function} in large
deviations theory, see Chapter 1.2 of \cite{Dem1998}. By Lemma 2.2.5
in \cite{Dem1998}, we have $$I_{n,\delta}\{(n-1)
E[\xi_{n,\delta}]\}=0.$$ Let $\overline{k}_{n,\delta}$ be
$(n-1)f_n(0)$ for the case of
$Var(\xi_{n,\delta})=E\xi_{n,\delta}^2-(E\xi_{n,\delta})^2=0$ (i.e.,
$\xi_{n,\delta}=f_n(0)$ is a degenerate random variable), and be a
solution of the equation $I_{n,\delta}(x)=\log n$ in $((n-1)
E[\xi_{n,\delta}],\infty)$ for the case of Var$(\xi_{n,\delta})>0$.
We will show that the solution
 uniquely exists in Section \ref{Sec_estimation}. For the simplicity of expression,
 we denote $\xi_n=\xi_{n,0}$,
$I_{n}=I_{n,0}$, and $\overline{k}_{n}=\overline{k}_{n,0}$.

\subsection{Notation}

In this paper, we investigate the flocking behavior of the system
(\ref{m2}) on the probability space $({\it \Omega},\mathcal{F},P)$,
where the sample space ${\it \Omega}=[0,1]^{dn}$, and the
$\sigma$-algebra $\mathcal{F}$ is the collection of all the Borel
subsets of ${\it \Omega}$. We say that a sequence of events $A_n\ (n
\geq 1)$ occurs with high probability (w.h.p.) if
$\lim_{n\rightarrow \infty} P(A_n)=1$.

A square matrix $M=(m_{ij})_{n\times n}$ is called stochastic, if
all elements $m_{ij}$ is nonnegative and for $1\leq i\leq n$,
$\sum_{j=1}^{n}m_{ij}=1$. For a matrix $A\in R^{n\times d}$, the
Frobenius norm $\|\cdot\|_F$ and the max norm $\|\cdot\|_{\max}$ of
the matrix $A$ are, respectively, defined as
$\|A\|_F=\sqrt{\sum_{i,j}a_{ij}^2}$ and
$\|A\|_{\max}:=\max_{i,j}|a_{ij}|$.

For two positive scalar sequences $g_1(n)$ and $g_2(n)$, we say that
(i) $g_1(n)=O(g_2(n))$ if there exists a constant $c>0$ and a
value $n_0>0$ such that $g_1(n)\leq c g_2(n)$ for any $n\geq n_0$;
(ii) $g_1(n)=\Theta(g_2(n))$ if there exist positive constants $c_1$ and
$c_2$ and a value $n_0>0$ such that $c_1g_2(n)\leq g_1(n)\leq c_2
g_2(n)$ for any $n\geq n_0$;
(iii) $g_1(n)=o(g_2(n))$ if $\lim_{n\rightarrow\infty}\frac{g_1(n)}{g_2(n)}=0$.

\addtocounter{footnote}{1}

\section{Main Results}\label{mainresult}
In this paper, we proceed with our analysis under the system \dref{m2} with the following
assumptions on the initial
states of all agents and the interaction function $f_n(\cdot)$:

 {A1) The initial positions $\{X_i(0)\}_{i=1}^n$ are independently
and uniformly distributed in $[0,1]^d$.}

 {A2) The nonlinear interaction weight $f_n(\cdot)$ is a
non-increasing integrable function satisfying{}\,(\ref{m0}), and the
interaction radius $r_n$ satisfies \dref{rad_con}.}

To avoid repetitive description we do not state the above assumptions in our results.

\subsection{Sufficient Conditions for Flocking}

Let $V_0:=\frac{1}{n}\sum_{i=1}^nV_i(0)$ and
$\overline{V}:=(V_0^{\rm T},V_0^{\rm T},\cdots,V_0^{\rm T})^{\rm
T}\in\mathbb{R}^{n\times d}$. Denote \bna \label{v0}
\mathcal{L}(V(0))=\|V(0)-\overline{V}\|_{\max}\left[\log\left(\frac{\|V(0)-\overline{V}\|_{F}}{\|V(0)-\overline{V}\|_{\max}}\right)+1\right].
\ena The sufficient condition for flocking can be stated as follows.

\begin{theorem}\label{result}
Suppose that the parameter $\alpha$ defined by $\dref{rad_con}$
satisfies $\alpha>\frac{2^{d-1}}{d\pi_d}$, and that there exist
positive constants $\delta$ and $\varepsilon$ such that
$\overline{k}_{n,\delta}\leq 1-\varepsilon$ holds for large $n$.
Then flocking is achieved w.h.p. if one of the following two
conditions holds:

{\rm (i)} $\alpha\varepsilon^d \leq (d+3)^{d/2}$, and \bna\label{v1}
\mathcal{L}(V(0))\leq c \delta n^2r_n \cdot \min \left\{r_n^{2d+2}
f_n^2((\delta+\varepsilon)r_n), \frac{
f_n^2\left(R_c+(\delta+\varepsilon)r_n \right)}{(\log
n)^{2d/(d-1)}}\right\}; \ena

{\rm (ii)} $\alpha\varepsilon^d > (d+3)^{d/2}$, and \bna \label{v2}
\mathcal{L}(V(0))\leq c \delta f_n^2[(\delta+\varepsilon)r_n]n^2
r_n\cdot \min\left\{r_n^{2d+2},1 \right\}, \ena where $c$ is a
positive constant depending on $d$, $\varepsilon$ and
$\alpha$\footnote{ Throughout this paper, a constant $c$ depending
on $\alpha$ means that $c$ depends on $\alpha$ only for
$\alpha<\infty$, and will not depend on $\alpha$ for
$\alpha=\infty$.}.
\end{theorem}

The proof of Theorem \ref{result} is put
in Subsection \ref{subsec_suff_proof}.

\begin{remark}
The value of \ $\overline{k}_{n,\delta}$ in Theorem \ref{result} can
be calculated by solving Equation{}\,(\ref{com4}), and some
theoretical results for $\overline{k}_{n,\delta}$  can also be
obtained, see Proposition \ref{Special} in
Section{}\,\ref{Sec_estimation}.
\end{remark}

\begin{remark} \label{remark32}It is clear that
$\mathcal{L}(V(0))$ defined via \dref{v0} satisfies $
\mathcal{L}(V(0))\leq \big(\frac{\log n}{2}+1\big)$ $ \max_{i,j}|v_{ij}(0)|.
$
Hence, Theorem \ref{result} still holds if $\mathcal{L}(V(0))$ in \dref{v1} and \dref{v2} is replaced by $\big(\frac{\log n}{2}+1\big)$ $\max_{i,j}|v_{ij}(0)|$, which means that under some conditions on the neighborhood radius and the interaction weights, the system (\ref{m2}) can reach flocking if the initial velocities are suitably small.
\end{remark}

In fact, under some further conditions on $f_n$, the parameter
$\delta$ in the condition of  Theorem{}\,\ref{result} can be
removed, see the following corollary.

\begin{corollary}\label{result2}
Let $\alpha>\frac{2^{d-1}}{d\pi_d}$ and $r_n=o(1)$. Suppose that there exists a
constant $\varepsilon>0$ such that  $ \overline{k}_{n} \leq
1-\varepsilon$ for large $n$, and
\begin{eqnarray}\label{S0_1}
\inf_{n}\left(\frac{1}{f_n(0)}\int_0^1 f_n(r_n y)y^{d-1}dy\right)=c_0>0.
\end{eqnarray}
Then the system $(\ref{m2})$ reaches flocking w.h.p., if one of the
following two conditions holds:

 {\rm (i)} $\alpha\varepsilon^d \leq
(d+3)^{d/2}$ and
\begin{eqnarray*}
\mathcal{L}(V(0))\leq c n^2r_n\cdot \min \left\{r_n^{2d+2} f_n^2\left(c_1\varepsilon r_n\right),
\frac{ f_n^2\left(R_c+c_1\varepsilon r_n
\right)}{(\log n)^{2d/(d-1)}}\right\};
\end{eqnarray*}

{\rm (ii)} $\alpha\varepsilon^d > (d+3)^{d/2}$ and
$\mathcal{L}(V(0))\leq c n^2r_n^{2d+3}f_n^2\left(c_1\varepsilon
r_n\right),$ where $c=c(d,\varepsilon,c_0,\alpha)$ and
$c_1=c_1(d,c_0,\alpha)$ are two positive constants.
\end{corollary}

The proof of Corollary \ref{result2} is put
in Subsection \ref{subsec_suff_proof}.

\begin{remark} From an intuitive point of view, the condition (\ref{S0_1}) means
that the interaction weight $f_n(\cdot)$ steadily decreases to zero, and cannot decay very fast.
\end{remark}

In the following, we provide an example to illustrate the result of Corollary \ref{result2}.

\begin{example}\label{example0}
Let $\beta>\frac{1}{d}$, $c'\in (0,1]$ and $\gamma>0$ be three constants. Set
\begin{eqnarray*}
 f_n(x):=\left\{%
\begin{array}{ll}
 c_n(1-x^{\gamma}/r_n^{\gamma}), \ \ & \mbox{if}~ x\leq r_n,\\
 0, &\mbox{otherwise},
\end{array}%
\right.
\end{eqnarray*}
where $r_n=n^{-1/d} \log^{\beta} n$, $c_n=\frac{c'}{\pi_d
\log^{d\beta} n}$. It is easy to verify that $f_n$ satisfies
(\ref{S0_1}). Recalling that $\xi_{n}=\xi_{n,0}=f_{n}(\|X-x_0\|),$
we obtain the following inequality:
\begin{eqnarray*}
\begin{aligned}
E\left[\xi_n\right]&=\int_{0}^{r_n}c_n \left(1-r_n^{-\gamma}x^{\gamma}\right)d\pi_d x^{d-1} dx=\frac{\gamma c_n \pi_d r_n^d}{\gamma+d}\leq \frac{\gamma }{(\gamma+d)n}.
\end{aligned}
\end{eqnarray*}
Using the following Proposition \ref{Special} (iii), we have for large $n$, $\overline{k}_n=nE[\xi_n](1+o(1))<1-\frac{d}{2(\gamma+d)}.$
Thus,
by Remark \ref{remark32} and Corollary \ref{result2} (ii), the system (\ref{m2}) reaches flocking w.h.p., if the velocity satisfies $\max_{i,j}|v_{ij}(0)|\leq c n^{\frac{-3}{d}} (\log n)^{3\beta-1}$ with $c>0$ being a positive constant.
\end{example}

\begin{figure}
\centering
\includegraphics[scale=0.43]{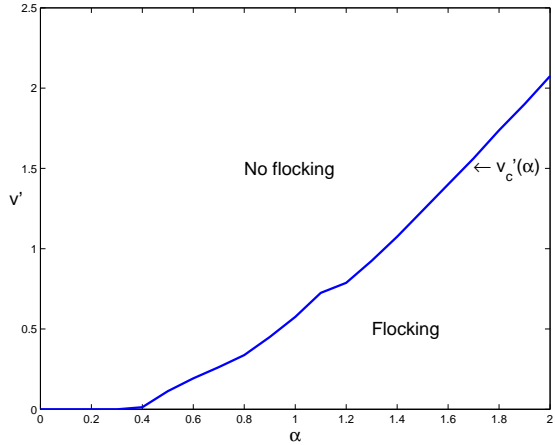}
\caption{There is a demarcation line $v_c'(\alpha)$ concerning with $v'$ and $\alpha$.  The flocking behavior can be reached below the line, while cannot be reached
above the line.}\label{Fig2}
\end{figure}

We also give some simulations for Theorem \ref{result} and  Corollary \ref{result2}. Assume the space dimension $d=2$, the interaction radius $r_n=\sqrt{\alpha \log n/n}$ with $\alpha$ being a positive constant, and
the weight function
\begin{eqnarray*}
 f_n(x)=\left\{%
\begin{array}{ll}
 \frac{1-x/r_n}{\alpha \pi\log n}, & \mbox{if}~x\leq r_n,\\
 0, &\mbox{otherwise}.
\end{array}%
\right.
\end{eqnarray*}
If $\alpha>1/\pi$, we can get
\begin{eqnarray}\label{simu1}
\begin{aligned}
 n^2r_n\cdot \min \left\{r_n^{6} f_n^2\left(\varepsilon r_n\right),
\frac{ f_n^2\left(\sqrt{\log n/(\pi n)}+\varepsilon r_n
\right)}{\log^4 n}\right\}=n^2 r_n^{7} f_n^2\left(\varepsilon r_n\right)=\Theta \left( \frac{\log n}{n} \right)^{3/2},
\end{aligned}
\end{eqnarray}
where $\varepsilon$ is a small positive constant.
The initial positions $\{X_i(0)\}_{i=1}^n$ are independently and uniformly chosen from the area $[0,1]^2$, and the initial velocities are set to be
\begin{eqnarray*}
 V_i(0)=\left\{%
\begin{array}{ll}
 (-v'n^{-\frac{3}{2}} \log^{\frac{1}{2}} n ,0), & \mbox{if}~X_{i1}(0)\leq 1/2\\
 (v'n^{-\frac{3}{2}} \log^{\frac{1}{2}} n ,0), & \mbox{otherwise}
\end{array}%
\right.,~~~~\forall 1\leq i\leq n,
\end{eqnarray*}
where $v'$ is a positive constant. By (\ref{v0}) we can obtain
\begin{eqnarray}\label{simu2}
\mathcal{L}(V(0))\approx v' n^{-\frac{3}{2}} \log^{\frac{1}{2}} n (\log n+1)=\Theta\left( n^{-\frac{3}{2}}  \log^{\frac{3}{2}} n  \right).
\end{eqnarray}
Simulations are carried out by choosing $n=600$, and the results are shown in Figure
\ref{Fig2}. It is shown that there is a demarcation line $v_c'(\alpha)$ concerning with $v'$ and $\alpha$ between the behaviors of flocking and no flocking. From (\ref{simu1}), (\ref{simu2}) and Figure \ref{Fig2}, the sufficient conditions for flocking in Theorem \ref{result} and  Corollary \ref{result2}
are tight in the order under this kind of interaction functions.


\subsection{Necessary Conditions for Flocking}

By  Remark \ref{remark32}, we see that for small initial velocity, the connectivity
of neighbor graphs can be maintained, and consequently the flocking behavior can be achieved.
 However, when the initial velocities become larger and larger, the connectivity of the dynamical
 neighbor graphs may be lost, even though the initial graph is connected. For such a case, the
 investigation for flocking becomes much harder. In order to simplify the analysis for the dependency
 of the flocking conditions on
initial velocities of the agents, we introduce the definition
 {$v$-flocking} as follows.
\begin{definition}\label{vconsensus}
Let $v>0$ be a constant. If  the system (\ref{m2}) reaches flocking
for any initial velocities satisfying $\max_{i\in
\mathcal{V}}\|V_i(0)\|\leq v$, then we say the system reaches
$v$-flocking.
\end{definition}

 From Definition \ref{vconsensus}, we see that the $v$-flocking behavior monotonically depends on $v$.
 For any initial positions of the agents, there exists a critical value $v_c$ such that the $v$-flocking
 can be achieved for $v<v_c$, and cannot be achieved for $v>v_c$. The investigation of $v_c$ is an
 interesting but very hard topic, and it falls into our future research.

The necessary conditions for $v$-flocking are presented as follows.

\begin{theorem}\label{result_0}
{\rm (i)}~If $\alpha<\frac{2^{d-1}}{d\pi_d}$,  then the system
$(\ref{m2})$ cannot achieve $v$-flocking w.h.p. for any $v>0$.

{\rm  (ii)}~If $\alpha \geq \frac{2^{d-1}}{d\pi_d}$ and $r_n=o(1)$,
then the system $(\ref{m2})$ cannot reach
$(2^{-d-1}\overline{k}_{n}r_n)$-flocking w.h.p. for
$\overline{k}_{n}={\it \Theta}(1)$ and $\overline{k}_{n}<2^d$, and
also cannot reach $(\frac{1}{2}\overline{k}_{n}r_n)$-flocking w.h.p.
for $\overline{k}_{n}=o(1)$.
\end{theorem}

The proof of Theorem \ref{result_0} is put in Subsection \ref{mainresultproof}.

 Under the definition of $v$-flocking, the sufficient conditions for flocking can be stated in a
 more clear manner.
\begin{corollary}\label{result_simple}
Assume that $\alpha>\frac{2^{d-1}}{d\pi_d}$ and $r_n=o(1)$, and that
there exists a constant $\varepsilon>0$ such that for large $n$,
{\rm (i)} $\varepsilon\leq \overline{k}_{n} \leq 1-\varepsilon$;
{\rm (ii)} $f_n(\varepsilon r_n)>\varepsilon f_n(0)$; and {\rm
(iii)}~$f_n\left(R_c+\varepsilon r_n\right)> r_n^2(\log
n)^{2d/(d-1)}n^{-2}$. Then the system $(\ref{m2})$ reaches $(c r_n^3
(\log n)^{-1})$-flocking w.h.p., where $c=c(d,\varepsilon,\alpha)>0$
is a constant.
\end{corollary}
The above corollary can be directly deduced from Corollary
\ref{result2}.

\begin{figure}
\centering
\includegraphics[scale=0.5]{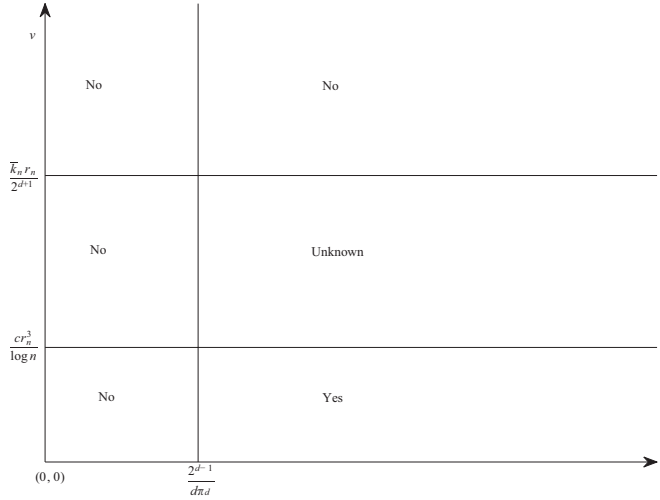}
\caption{``Yes" (or ``No") means that in a probability sense, the
system (\ref{m2}) can (cannot) reach $v$-flocking, and ``unknown"
means that we cannot judge whether the system
 reaches $v$-flocking}\label{Fig1}
\end{figure}

From Theorem \ref{result_0} (i) and Corollary \ref{result_simple},
we see that in a probability sense, $R_c$ can be considered as the
smallest possible interaction radius for flocking, and also the
critical interaction radius for $v$-flocking with $v\leq c r_n^3
(\log n)^{-1}$. We illustrate the result of  Theorem \ref{result_0}
and Corollary \ref{result_simple} in Figure \ref{Fig1} where the
interaction function $f_n(\cdot)$ satisfies the conditions of
Corollary{}\,\ref{result_simple}.

\section{Calculation and Estimation of $\overline{k}_{n,\delta}$ }\label{Sec_estimation}
In this section, we present the approach to the estimation and
calculation of $\overline{k}_{n,\delta}$. We only consider the case
of Var$(\xi_{n,\delta})>0.$ Note that $I_{n,\delta}((n-1)
E[\xi_{n,\delta}])=0$, and
\begin{eqnarray*}\label{com2}
\lim_{x\rightarrow\infty} I_{n,\delta}(x)\geq
\lim_{x\rightarrow\infty}\left\{ x-(n-1) \log \left(E\left[{\rm
e}^{\xi_{n,\delta} } \right]\right)\right\}=\infty.
\end{eqnarray*}
By the continuity of $I_{n,\delta}(x)$, we see that the solution of
the equation $I_{n,\delta}(x)=\log n$ in $((n-1)
E[\xi_{n,\delta}],\infty)$ exists. Furthermore, by Cauchy-Schwarz
inequality, we have
\begin{eqnarray}\label{com2_add}
\begin{aligned}
&\frac{d^2}{d^2 \theta}\log \left(E\left[{\rm
e}^{\theta\xi_{n,\delta} } \right]\right)=\frac{E\left[{\rm
e}^{\theta\xi_{n,\delta} } \right]E\left[\xi_{n,\delta}^2{\rm
e}^{\theta\xi_{n,\delta} } \right]-\left(E\left[\xi_{n,\delta} {\rm
e}^{\theta\xi_{n,\delta} } \right] \right)^2 }{\left(E\left[{\rm
e}^{\theta\xi_{n,\delta} } \right] \right)^2}>0.
\end{aligned}
\end{eqnarray}

\noindent Hence, the following equation
\begin{eqnarray}\label{com4_add}
\begin{aligned}
\frac{d}{d\theta}\left\{\theta x-(n-1) \log \left(E\left[{\rm
e}^{\theta\xi_{n,\delta} }
\right]\right)\right\}=x-\frac{(n-1)E\left[\xi_{n,\delta} {\rm
e}^{\theta\xi_{n,\delta} } \right]}{E\left[{\rm
e}^{\theta\xi_{n,\delta} } \right]}=0
\end{aligned}
\end{eqnarray}
has  a unique solution $\theta^*(x)$ for any $x>0$. Moreover, by (\ref{com2_add}), the right-hand side of \dref{com4_add2}
is strictly monotonically increasing. Thus,
the equation
\begin{eqnarray}\label{com4_add2}
\ay\begin{array}[b]{lll} \log n=I_{n,\delta}(x)&=&\theta^*(x)
x-(n-1) \log
\left(E\left[{\rm e}^{\theta^*(x)\xi_{n,\delta} } \right]\right)\\
&=&\displaystyle \theta^*(x) \frac{(n-1)E\left[\xi_{n,\delta} {\rm
e}^{\theta^*(x) \xi_{n,\delta} } \right]}{E\left[{\rm
e}^{\theta^*(x) \xi_{n,\delta} } \right]}-(n-1) \log
\left(E\left[{\rm e}^{\theta^*(x)\xi_{n,\delta} } \right]\right)
\end{array}
\end{eqnarray}
with respect to $\theta^*(x)$ also has a unique solution $\overline{\theta}_{n,\delta}$. By (\ref{com4_add}), we know that $\overline{k}_{n,\delta}$ is the unique solution of $I_{n,\delta}(x)=\log n$ in $((n-1)
E[\xi_{n,\delta}],\infty)$. Combining (\ref{com4_add2}) with (\ref{com4_add}), we see that $\overline{k}_{n,\delta}$ and $\overline{\theta}_{n,\delta}$ can be obtained by solving the following equations,
\begin{eqnarray}\label{com4}
\begin{aligned}
&\frac{(n-1)E\left[\xi_{n,\delta} {\rm
e}^{\overline{\theta}_{n,\delta}\xi_{n,\delta} }
\right]}{E\left[{\rm e}^{\overline{\theta}_{n,\delta}\xi_{n,\delta}
} \right]}=\frac{\log n + (n-1) \log \left(E\left[{\rm
e}^{\overline{\theta}_{n,\delta}\xi_{n,\delta} }
\right]\right)}{\overline{\theta}_{n,\delta}}=\overline{k}_{n,\delta}.
\end{aligned}
\end{eqnarray}

Furthermore, if $f_n$ does not decay very fast, there are some theoretical results to estimate
$\overline{k}_{n,\delta}$, $\overline{k}_{n}(=: \overline{k}_{n,0})$ and $
\overline{\theta}_{n}(=: \overline{\theta}_{n,0})$.

\begin{proposition}\label{Special}
Suppose that $r_n=o(1)$, and $f_n$ satisfies $(\ref{S0_1})$. Let
$c_0$ be the constant appearing in $(\ref{S0_1})$. Then the
following results hold:

{\rm (i)} $\overline{k}_n={\it \Theta}(nr_n^d f_n(0))$ and
$\overline{\theta}_n=O(1/f_n(0));$

{\rm (ii)} There exists a constant $c_1=c_1(d,\alpha,c_0)>0$, such
that for any $\delta\in(0,1)$, $\overline{k}_{n,\delta}\leq
(1+c_1\delta) \overline{k}_n;$

{\rm  (iii)} For the case of $\alpha=\infty$,
$\overline{k}_n=nE[\xi_n](1+o(1))$ and
$\overline{\theta}_n=o(1/f_n(0)).$
\end{proposition}

The proof of Proposition \ref{Special} is in Appendix \ref{probosition}.

Using a similar method as that of Proposition \ref{Special}, we obtain the following results.

\begin{corollary}\label{remark_eta}
Suppose $r_n=o(1)$. Then for any constant $\delta>0$, we have
$\overline{k}_{n,\delta}={\it \Theta}(nr_n^d f_n(0))$ and
$\overline{\theta}_{n,\delta}=O(1/f_n(0))$. If $\alpha=\infty$, then
$\overline{k}_{n,\delta}=nE[\xi_{n,\delta}](1+o(1))$ and
$\overline{\theta}_{n,\delta}=o(1/f_n(0))$.
\end{corollary}

In the following, we present two examples to show how to calculate the value of
$\overline{k}_n$. For this, we introduce some notations in large deviations theory.  Define
$H:[0,\infty)\rightarrow \mathbb{R}$ by $H(0)=1$ and
\begin{eqnarray*}
H(a)=1-a+a\log a,~~~~~a>0.
\end{eqnarray*}
Note that $H(1)=0$ and the unique turning point of $H$ is the
minimum at $1$. Also $H(a)/a$ is increasing on $(1,\infty)$. Let
$H_-^{-1}:[0,1]\rightarrow[0,1]$ be the unique inverse of the
restriction of $H$ to $[0,1]$, and let
$H_+^{-1}:[0,\infty)\rightarrow[1,\infty)$ be the inverse of the
restriction of $H$ to $[1,\infty)$.

\begin{example}\label{example1}
Set $r_n=c n^{-1/d} (\log n)^{1/d}$ and
 \begin{eqnarray*}
 f_n(x):=\left\{%
\begin{array}{ll}
 b_n, & x\leq r_n,\\
 0, &x> r_n,
\end{array}%
\right.
\end{eqnarray*}
where $c$ and $b_n$ are positive constants. Recall that
$\xi_{n}=\xi_{n,0}=f_{n}(\|X-x_0\|),$ where
$x_0=(\frac{1}{2},\frac{1}{2},\cdots,\frac{1}{2})\in\mathbb{R}^d$
and $X$ is a random variable distributed uniformly in $[0,1]^d$. To
solve (\ref{com4}) we first calculate the following equation,
\begin{eqnarray}\label{ex1_1}
\ay\begin{array}[b]{lll} E\left[{\rm e}^{\overline{\theta}_n \xi_n}
\right]&=& P(\|X-x_0\|\leq r_n) {\rm e}^{\overline{\theta}_n b_n} + P(\|X-x_0\|> r_n)\\
&=& 1+\left({\rm e}^{\overline{\theta}_n b_n}-1 \right)\pi_d c^d
n^{-1} \log n.
\end{array}
\end{eqnarray}
Similarly,
\begin{eqnarray}\label{ex1_2}
E\left[\xi_n {\rm e}^{\overline{\theta}_n \xi_n} \right]={\rm
e}^{\overline{\theta}_n b_n}\pi_d c^d b_n n^{-1}\log n.
\end{eqnarray}
By Proposition \ref{Special} (i), we have $\overline{\theta}_n
b_n=\overline{\theta}_n f_n(0)=O(1)$. Then by (\ref{ex1_1}), we
obtain $E[{\rm e}^{\overline{\theta}_n \xi_n}]=1+o(n^{-1/2})$ and
\bna \label{ex1_3} \log E\left[{\rm e}^{\overline{\theta}_n \xi_n}
\right]=\left({\rm e}^{\overline{\theta}_n b_n}-1 \right)\pi_d c^d
n^{-1} \log n+O(n^{-2}\log^2 n). \ena Substituting (\ref{ex1_2}) and
\dref{ex1_3} into (\ref{com4}), we have
\begin{eqnarray}\label{ex1_4}
\overline{\theta}_n b_n {\rm e}^{\overline{\theta}_n b_n}- {\rm
e}^{\overline{\theta}_n b_n}=\left( \frac{1}{\pi_d c^d}-1
\right)\left(1+o\left(n^{-1/2}\right)\right)+o\left(n^{-1/2}\right).
\end{eqnarray}
By the definition of $H_+^{-1}$, we see that
$\log(H_+^{-1}(\frac{1}{\pi_d c^d } ))$ is the unique solution of
the equation $x {\rm e}^{x}- {\rm e}^{x}= \frac{1}{\pi_d c^d}-1$
with respect to $x$. Substituting this into (\ref{ex1_4}) we have
$\overline{\theta}_n=\frac{1}{b_n}\log(H_+^{-1}(\frac{1}{\pi_d c^d }
))(1+o(1))$. Furthermore, by (\ref{com4}) we have
$\overline{k}_n=H_+^{-1}(\frac{1}{\pi_d c^d })\pi_d c^d b_n \log n
(1+o(1)).$
\end{example}

\begin{example}
Let $r_n=c n^{-1/d} (\log n)^{1/d}$ and
 \begin{eqnarray*}
 f_n(x):=\left\{%
\begin{array}{ll}
 b_n\left(1-r_n^{-1}x\right), & x\leq r_n,\\
 0, &x>r_n,
\end{array}%
\right.
\end{eqnarray*}
where  $c$ and $b_n$ are positive constants.
Similar to (\ref{ex1_1}) and (\ref{ex1_2}) we have
\begin{eqnarray*}
\begin{aligned}
&E\left[{\rm e}^{\overline{\theta}_n \xi_n}
\right]=1-\pi_d r_n^d+\int_{0}^{r_n} {\rm e}^{\overline{\theta}_n b_n \left(1-r_n^{-1}x\right)}d\pi_d x^{d-1} dx\\
&\hspace{1.5cm}=1-\pi_d r_n^d+d\pi_d r_n^d \int_0^{1}{\rm e}^{\overline{\theta}_n b_n y}(1-y)^{d-1}dy ,\\
&E\left[\xi_n {\rm e}^{\overline{\theta}_n \xi_n}
\right]=\int_{0}^{r_n}b_n \left(1-r_n^{-1}x\right)
{\rm e}^{\overline{\theta}_n b_n (1-r_n^{-1}x)}d\pi_d x^{d-1} dx\\
&\hspace{1.8cm}=d\pi_d r_n^d b_n \int_0^{1}y{\rm
e}^{\overline{\theta}_n b_n y}(1-y)^{d-1}dy .
\end{aligned}
\end{eqnarray*}
By Proposition \ref{Special} (i), we have $\overline{\theta}_n
b_n=\overline{\theta}_n f_n(0)=O(1)$. Similar to (\ref{ex1_4}) we
have \ay\begin{eqnarray}\label{ex2_1} &&\overline{\theta}_n b_n
\int_0^{1}y{\rm e}^{\overline{\theta}_n b_n
y}(1-y)^{d-1}dy-\int_0^{1}{\rm e}^{\overline{\theta}_n b_n y}(1-y)^{d-1}dy\nonumber\\
&=& \left(\frac{1}{d\pi_d
c^d}-\frac{1}{d}\right)\left(1+o\left(n^{-1/2}\right)\right)+o\left(n^{-1/2}\right).
\end{eqnarray}
Let $x^*$ be the unique solution of the equation
\begin{eqnarray*}\label{ex2_2}
\begin{aligned}
x  \int_0^{1}y{\rm e}^{x y}(1-y)^{d-1}dy-\int_0^{1}{\rm e}^{x
y}(1-y)^{d-1}dy= \frac{1}{d\pi_d c^d}-\frac{1}{d}
\end{aligned}
\end{eqnarray*}
with respect to $x$. Substituting $x^*$ into (\ref{ex2_1}) we can get
$\overline{\theta}_n = \frac{x^*}{b_n}(1+o(1))$. Substituting this into
 (\ref{com4}) we
obtain
$$\overline{k}_n=d\pi_d c^d b_n \log n\left(\int_0^{1}y {\rm e}^{x^*
y}(1-y)^{d-1}dy\right)
 (1+o(1)).$$
\end{example}

\section{Proofs of Main Results}\label{Sec_proof_mrs}

\subsection{Key Lemmas}\label{Sec_Interact}

In order to prove the main results of this paper, we first present
some key lemmas. By{}\,(\ref{m2}) we have $V(t+1)=P(t)\cdots
P(0)V(0)$. To study the flocking behavior, we need to deal with the
convergence of the matrix product $P(t)\cdots P(0)$. Compared with
our previous work \cite{Chen2014}, a crucial point lies in the fact
that the weighted average matrix $P(t)$ is determined by the
nonlinear function $f_n(\cdot)$ and the distance between agents. We
apply the large deviations techniques to deal with this.

Introduce the maximum weighted degree of all agents at time $t$ as
$$\Delta_{n}(t):=\max_{1\leq i \leq n}\Big\{\sum_{1\leq j \leq n, j\neq i}f_{n}(\|X_i(t)-X_j(t)\|)\Big\}.$$  Note that  $f_n(\cdot)$ is a non-increasing function. If there exists a positive constant $\delta$ such that
\begin{eqnarray}\label{newsec_2}
\|X_i(t)-X_j(t)\| \geq \|X_i(0)-X_j(0)\|-\delta r_n,
\end{eqnarray}
then by the definition of $f_{n,\delta}(\cdot)$ we have
\begin{eqnarray}\label{newsec_3}
f_{n}(\|X_i(t)-X_j(t)\|)\leq f_{n,\delta}(\|X_i(0)-X_j(0)\|).
\end{eqnarray}
Thus,
we have the following lemma of $\Delta_{n}(t)$.

\begin{lemma}\label{lemma_lzx}
Assume that there exists a positive constant $\delta$, such that
$$\Delta_{n,\delta}:=\max_{1\leq i \leq n}\left\{\sum_{1\leq j \leq n, j\neq i}f_{n,\delta}(\|X_i(0)-X_j(0)\|)\right\}\leq 1,$$ and that
the inequality $(\ref{newsec_2})$ holds for all $i$ and $j$ and all
$t\geq 0$, then we have $\Delta_{n}(t)\leq 1$.
\end{lemma}

It is clear that under the conditions of Lemma \ref{lemma_lzx},
$P(t)\ (t\geq 0)$ is a stochastic matrix. The stochastic matrices
may bring convenience for the convergence property of $P(t)\cdots
P(0)$.


For $\Delta_{n,\delta}$, we have the following results.
\begin{lemma}\label{Delta1}
Given a constant $\delta\geq 0$, suppose that ${\rm
Var}(\xi_{n,\delta})>0$ holds. Then,
\begin{eqnarray}\label{D1_00}
\Delta_{n,\delta} \leq
\overline{k}_{n,\delta}(1+o(1)),~~~~\mbox{w.h.p.}
\end{eqnarray}
\end{lemma}

The large deviations techniques are used in the proof of Lemma \ref{Delta1}, and the proof details are put in Appendix \ref{Delta1proof}.

In the following, we study the eigenvalues of the matrix $P(t)$.
Under the conditions of Lemma \ref{lemma_lzx}, $P(t) (t\geq 0)$ is a
symmetric stochastic matrix. Thus, the eigenvalues of $P(t)$ are all
real numbers. We denote $\lambda_i(t)=\lambda_i(t,n)$ as the
$i$-largest eigenvalues of $P(t)$, and arrange all eigenvalues
according to the following order,
$$1=\lambda_1(t)\geq\lambda_2(t)\geq
\cdots\geq\lambda_n(t)\geq-1.$$ The essential spectral radius
$\bar{\lambda}(t)$ of $P(t)$ is defined as
$\bar{\lambda}(t)=\bar{\lambda}(t,n):=\max\{|\lambda_2(t)|,$
$|\lambda_n(t)|\}$, which plays a key role for the convergence of
the matrix product $P(t)\cdots P(0)$.

\begin{lemma}\label{Eigenvalue}
Assume that $\alpha$ defined by $(\ref{rad_con})$ satisfies
$\alpha>\frac{2^{d-1}}{d\pi_d}$, and that there exists a constant
$\delta>0$, such that w.h.p. the following inequality holds,
\begin{eqnarray}\label{Ei_0}
\sup_{i,j,t}\left\{\|X_i(t)-X_j(t)-X_i(0)+X_j(0)\|\right\} \leq
\delta r_n.
\end{eqnarray}
Furthermore, for any given constant $\varepsilon >0$,
$\overline{k}_{n,\delta}$ satisfies $ \overline{k}_{n,\delta} \leq
1-\varepsilon$. Then the essential spectral radius
$\bar{\lambda}(t)\ (t>0)$ satisfies the following inequalities:

{\rm (i)}~ If $\alpha\varepsilon^d \leq (d+3)^{d/2}$, then w.h.p.
\begin{eqnarray*}
\overline{\lambda}(t)\leq 1-cn^2 \cdot \min \left\{r_n^{2d+2}
f_n^2((\delta+\varepsilon)r_n), \frac{
f_n^2\left(R_c+(\delta+\varepsilon)r_n \right)}{(\log n)^{2d/(d-1)}}
\right\};
\end{eqnarray*}

{\rm (ii)}~ If $\alpha\varepsilon^d > (d+3)^{d/2}$, then w.h.p.
\begin{eqnarray*}
\overline{\lambda}(t)\leq 1-cn^2
f_n^2((\delta+\varepsilon)r_n)\min\left\{r_n^{2d+2},1 \right\},
\end{eqnarray*}
where $c$ is a positive constant depending on $\alpha$, $\varepsilon$ and $d$.
\end{lemma}

The proof of Lemma \ref{Eigenvalue} is put in Appendix
\ref{eigenvalueproof}.

\vskip 3mm

\subsection{Proofs of  Theorem \ref{result} and Corollary \ref{result2}}\label{subsec_suff_proof}

\emph{Proof~of Theorem} \ref{result}~ By the fact that the row sum
of $P(t)$ is $1$, we  see that  $P(t)\overline{V}=\overline{V}$. By
(\ref{m2}), we have
\begin{eqnarray}\label{r1}
\begin{aligned}
V(k+1)&=P(k)P(k-1)\cdots P(0)V(0)=P(k)P(k-1)\cdots P(0)(V(0)-\overline{V})+\overline{V}.
\end{aligned}
\end{eqnarray}
Using (\ref{r1}) and Corollary 1 in \cite{Smale2007}, we obtain that
\begin{eqnarray}\label{r2}
\|V(k)-\overline{V}\|_F \leq
\|V(0)-\overline{V}\|_F\prod_{i=0}^{k-1} \overline{\lambda}(i).
\end{eqnarray}

Denote
\begin{eqnarray*}
\lambda^*=\left\{%
\begin{array}{ll}
 1-cn^2 \min \left\{r_n^{2d+2}
f_n^2((\delta+\varepsilon)r_n), \dfrac{
f_n^2\left(R_c+(\delta+\varepsilon)r_n \right)}{(\log n)^{2d/(d-1)}}
\right\}, &\mbox{if  $\alpha\varepsilon^d \leq (d+3)^{d/2}$},\\
 1-cn^2
f_n^2((\delta+\varepsilon)r_n)\min\left\{r_n^{2d+2},1
\right\}, &\mbox{else},
\end{array}%
\right.
\end{eqnarray*}
where $c$ is the same constant appearing in Theorem
\ref{Eigenvalue}. Take
 $$k^*=\left\lceil
\log_{\frac{1}{\lambda^*}}\left(\frac{
\|V(0)-\overline{V}\|_{F}}{\|V(0)-\overline{V}\|_{\max}} \right)
\right\rceil.
$$ It is easy to see that $k^*$ satisfies $k^*\leq \frac{2}{1-\lambda^*}
\log\left(\frac{
\|V(0)-\overline{V}\|_{F}}{\|V(0)-\overline{V}\|_{\max}} \right)$.
By the condition of the theorem, we have
\begin{eqnarray}\label{r9}
\|V(0)-\overline{V}\|_{\max}\left(\log\left(\frac{\|V(0)-\overline{V}\|_{F}}{\|V(0)-\overline{V}\|_{\max}}\right)
+1 \right)\leq \frac{(1-\lambda^*)\delta r_n}{4\sqrt{d}},
\end{eqnarray}
then
\begin{eqnarray}\label{r4}
\begin{aligned}
2\|V(0)-\overline{V}\|_{\max}\left(\sqrt{d}k^*+
\frac{1}{1-\lambda^*} \right)\leq \delta r_n.
\end{aligned}
\end{eqnarray} We assert that if (\ref{r4}) holds, then  we have w.h.p.
\begin{eqnarray}\label{r5}
\max_{1\leq i,j \leq n}d_{ij}(k) \leq \delta r_n,~~~~~~~~\forall
k\geq 0,
\end{eqnarray}
where $d_{ij}(k)=\|X_i(k)-X_j(k)-X_i(0)+X_j(0)\|.$

We prove the aforementioned assertion by induction.
It is clear that \dref{r5} holds for $k=0$. We assume that the assertion holds for $0<s\leq k$, that is, w.h.p.
\begin{eqnarray*}\label{r7}
\max_{i,j} d_{ij}(s)\leq \delta r_n,~~~~~~\forall\ 0\leq s\leq k.
\end{eqnarray*}
By Lemmas \ref{lemma_lzx} and \ref{Delta1}, the matrix $P(s)\ (0\leq
s\leq k)$ is stochastic. Denote
$$a(t):=\Bigg(\sum_{j=1}^d\left(\max_{1\leq i_1,i_2 \leq n}|v_{i_1,
j}(t)-v_{i_2,j}(t)|^2 \right) \Bigg)^{1/2}.$$ It is clear that the
sequence $a(t)$ is monotonically non-increasing, and thus
\begin{eqnarray}\label{r3}
a(t)\leq a(0) \leq 2\sqrt{d} \|V(0)-\overline{V}\|_{\max}.
\end{eqnarray} Moreover, using Lemma \ref{Eigenvalue}, we have $\overline{\lambda}(s) \leq \lambda^*$ for $s\leq
k$. By (\ref{r2}) and (\ref{r3}), we obtain that w.h.p.
\ay\bna\label{r8} \max_{i,j}\left\{d_{ij}(k+1)\right\}&=&
\max_{i,j}\left\{\big\|\sum_{s=0}^{k-1} \left(
V_i(s)-V_j(s)\right)\big\|\right\}\nn\\ &\leq&
k^*a(0)+2\sum_{s=k^*}^k \left(\max_{i}\big\|V_i(s)-V_0\big\|\right)
\leq k^*a(0)+2\sum_{s=k^*}^k
\|V(s)-\overline{V}\|_F  \nn\\
 &\leq& k^*2\sqrt{d} \|V(0)-\overline{V}\|_{\max}+2 \|V(0)-\overline{V}\|_F\sum_{s=k^*}^k
(\lambda^*)^s\nn\\
&\leq& 2 \|V(0)-\overline{V}\|_{\max}\left(\sqrt{d}k^*+
\frac{1}{1-\lambda^*} \right)\leq \delta r_n. \ena The assertion
holds for $s=k+1$. By induction, the assertion holds for all $s\geq
0$.

By (\ref{r2}), (\ref{r5}) and Lemma \ref{Eigenvalue}, w.h.p.
\begin{eqnarray*}
\lim_{k\rightarrow\infty}\|(V(k)-\overline{V})\|_F \leq
\|V(0)-\overline{V}\|_F \cdot \lim_{k\rightarrow\infty}
(\lambda^*)^k=0
\end{eqnarray*}
holds for all large $n$. Substituting the value of $\lambda^*$ into
(\ref{r9}) yields our result.   \rulex

\emph{Proof of Corollary} \ref{result2}~ By Proposition
\ref{Special}, we see that there exists a constant
$c'=c'(d,\alpha,c_0)>0$, such that for any $\delta\in(0,1)$, we have
$\overline{k}_{n,\delta}\leq (1+c'\delta) \overline{k}_n$. Take
$\delta=\frac{\varepsilon}{2c'}$,  we have $ \overline{k}_{n,\delta}
\leq \left(1+\frac{\varepsilon}{2}\right) \overline{k}_n
<1-\frac{\varepsilon}{2}.$ Let $c_1=1+\frac{1}{2c'}$, the result
follows by the same argument as that of Theorem \ref{result}. \rulex

\subsection{Proof of Theorem \ref{result_0}}\label{mainresultproof}
(i) Set $M_n:=\lfloor 1/(2R_c) \rfloor-1$, where
$R_c=\sqrt[d]{\frac{2^{d-1}\log n}{d \pi_d n}}$, and  $\lfloor x
\rfloor$ denotes the largest integer no bigger than $x$. For any
integer $k\in[0,M_n]$, define the point
$x_k:=((2k+1)R_c,0,\cdots,0)\in [0,1]^d$ and set the event
$$A_k:=\left\{\mathcal{X}_n\cap B(x_k;n^{-\frac{1}{d}})\neq \emptyset,  \mathcal{X}_n\cap\left[B(x_k; r_n+n^{-\frac{1}{d}})\setminus B(x_k;n^{-\frac{1}{d}})\right]=\emptyset  \right\}.$$
Let $B(x;r)$ denotes the $d$-dimensional ball centered at $x$ with
the radius $r$. By the similar analysis as that of the equation
(4.13) in \cite{Chen2014}, we have for large $n$
\begin{eqnarray*}
\begin{aligned}
P\left(\bigcup_{0\leq k\leq M_n}A_k\right)&>1-2{\rm e}^{-n^{1/4}}-\left[1-\left(1-\exp\left\{-\frac{(n-n^{\frac{3}{4}})\pi_d}{2^{d-1}n}\right\} \right)\right.\\
&~~~~~~~~~~~~~~ \ \ \ ~~\ \ \ \cdot\left.\exp\left\{-\frac{(n+n^{\frac{3}{4}})\pi_d[(r_n+n^{-\frac{1}{d}})^d-\frac{1}{n}]}{2^{d-1}}\right\}\right]^{M_n+1}\\
&>1-2{\rm e}^{-n^{1/4}}-\left[1-\frac{\pi_d}{2^{d}} \exp\left\{-\frac{\pi_d\alpha'\log n}{2^{d-1}}\right\}\right]^{M_n+1}\\
&>1-2{\rm e}^{-n^{1/4}}-\exp\left\{-n^{-\frac{\pi_d
\alpha''}{2^{d-1}}+\frac{1}{d}}   \right\},
\end{aligned}
\end{eqnarray*}
where the constants $\alpha',\alpha''$ satisfy
$\alpha<\alpha'<\alpha''<\frac{2^{d-1}}{d\pi_d}$. So, w.h.p. there
must exist an integer $k\in [0, M_n]$ such that $A_k$ happens.
Without loss of generality we assume $A_0$ happens, then there exist
some agents lying in $B(x_0;n^{-1/d})$ which have no neighbor in
$[0,1]^d\setminus B(x_0;n^{-1/d})$ at the initial time. Taking their
initial velocities to be $(0,-v,\cdots,-v)$ and the other agents'
initial velocities to be $(0,v,\cdots,v)$ for any $v>0$, we see that
the system (\ref{m2}) cannot achieve flocking w.h.p.

(ii)~ Let $\widetilde{i}$ be the index of the agent whose position
satisfies $\|X_{\widetilde{i}}(0)\|\leq \|X_{i}(0)\|$ for any $
1\leq i \leq n.$ Then, for any constant $\varepsilon>0$ we have
\begin{eqnarray*}
\begin{aligned}
P\left(\|X_{\widetilde{i}}(0)\|>2\bigg(\frac{\varepsilon\log n}{\pi_d n}\bigg)^{1/d} \right)&
=P\left(\bigcap_{i=1}^n\left\{\|X_{i}(0)\|>2\bigg(\frac{\varepsilon\log n}{\pi_d n}\bigg)^{1/d}\right\}  \right)\\
&=\left(1-\frac{\varepsilon\log n}{n}\right)^n=n^{-\varepsilon}(1+o(1)),
\end{aligned}
\end{eqnarray*}
which indicates that $\|X_{\widetilde{i}}(0)\| \leq
2(\frac{\varepsilon\log n}{\pi_d n})^{1/d}$ holds w.h.p. For some
constant $v_0$, we set
$V_{\widetilde{i}}(0)=(-\frac{v_0}{\sqrt{d}},\cdots,-\frac{v_0}{\sqrt{d}})^{\rm
T}$, and $V_{i}(0)=-V_{\widetilde{i}}(0)$ for $i\neq \widetilde{i}$.

If $\overline{k}_n={\it \Theta}(1)$, we take
$v_0=2^{-d-1}\overline{k}_n r_n$. By (\ref{m2}), w.h.p. all elements
of $X_{\widetilde{i}}(1)=X_{\widetilde{i}}(0)+V_{\widetilde{i}}(0)$
are less than $0$, and $\|X_{\widetilde{i}}(1)-X_{i}(1)\|\geq
2^{-d}\overline{k}_n r_n$ holds for $i\neq \widetilde{i}$. By Lemma
\ref{Delta1}, it is easy to see that there exists a small constant
$\epsilon>0$ such that
$$d_{\widetilde{i}}(1)=\sum_{i\neq \widetilde{i}}f_n(\|X_{\widetilde{i}}(1)-X_i(1)\|)\leq \frac{(1-\epsilon)\overline{k}_n}{2^d}~~~~\mbox{w.h.p.}.$$
Thus, for $i\neq \widetilde{i}$, \ay\begin{eqnarray*}
&&\|X_{\widetilde{i}}(2)-X_{i}(2)\|\\
&=&\|X_{\widetilde{i}}(1)+[1-2d_{\widetilde{i}}(1)]V_{\widetilde{i}}(0) -X_{i}(1)+[1-2f_n(\|X_{\widetilde{i}}(1)-X_{i}(1)\|)]V_{\widetilde{i}}(0)\|\\
&\geq &2^{-d}\overline{k}_n r_n \left(1+\left[1-
\bigg(1-\frac{\epsilon}{2}\bigg) \frac{\overline{k}_n}{2^d}\right]
\right) \quad \mbox{w.h.p.}
\end{eqnarray*}
Repeating the above process we have for all $t>1$,
\begin{eqnarray*}
\begin{aligned}
&\|X_{\widetilde{i}}(t)-X_{i}(t)\|\\
\geq& 2^{-d}\overline{k}_n r_n \left(1+ \left[1-
\bigg(1-\frac{\epsilon}{2}\bigg)
\frac{\overline{k}_n}{2^d}\right]+\cdots+ \left[1-
\bigg(1-\frac{\epsilon}{2}\bigg)
\frac{\overline{k}_n}{2^d}\right]^{t-1}\right)\quad  \mbox{w.h.p.}
\end{aligned}
\end{eqnarray*}
 Therefore, there exists a time instant $T$ such that w.h.p. for all $i\neq \widetilde{i}$: (i)
$\|X_{\widetilde{i}}(T)-X_{i}(T)\|\geq r_n$; (ii) all elements of
$V_{\widetilde{i}}(T)$ are negative; (iii) all elements of
$V_{i}(T)$ are positive. For such a case, the system cannot reach
flocking.

If $\overline{k}_n=o(1)$, we take $v_0=\frac{1}{2}\overline{k}_n
r_n$. Similarly we see that there exists a constant $\epsilon>0$
such that w.h.p. $d_{\widetilde{i}}(1)\leq
(1-\epsilon)\overline{k}_n$, and
\begin{eqnarray*}
\begin{aligned}
\|X_{\widetilde{i}}(t)-X_{i}(t)\|\geq \overline{k}_n r_n \left(1+
\left[1-
\bigg(1-\frac{\epsilon}{2}\bigg)\overline{k}_n\right]+\cdots+
 \left[1- \bigg(1-\frac{\epsilon}{2}\bigg) \overline{k}_n\right]^{t-1}\right)
\end{aligned}
\end{eqnarray*}
holds for $t>1$. Thus, the system cannot reach flocking.

\section{Concluding Remarks}\label{conclusion}
A fundamental issue in the investigation of multi-agent systems is how the local interactions affect the collective behavior of the overall systems.
This paper studied a  discrete-time nonlinear multi-agent system, where
the nonlinear interaction function decays according to the distance between agents. By applying large deviations techniques to estimate the
essential spectral gap of average matrices whose elements are determined by the nonlinear function and the distance between agents, we provide sufficient conditions and necessary conditions for the flocking behavior.
Some interesting problems deserve to be further investigated, for example, how to obtain the
critical value of $v$ for $v$-flocking of our multi-agent model, and how to analyze the flocking behavior of the continuous-time multi-agent models.

\section*{Appendix}


\section*{Proof of Lemma \ref{rgg}}\label{App_A}

The lemma holds for $d=2$, see \cite{kumar1998} and
\cite{penrose1997}. In the following, we prove the lemma for $d\geq
3$.

For $\alpha>\frac{2^{d-1}}{d\pi_d}$, by Theorem 7.14 in
\cite{penrose2003}, we see that the minimum degree of
$G(\mathcal{X}_n;r_n)$ is equal to ${\it \Theta}(nr_n^d)$ w.h.p.,
which indicates that there is no isolated vertex in
$G(\mathcal{X}_n;r_n)$ w.h.p. Following the proof idea of Theorem
3.1 in \cite{kumar1998}, we see that for $d\geq 3$, the probability
that $G(\mathcal{X}_n;r_n)$ is not connected has the same order as
the probability that the graph $G(\mathcal{X}_n;r_n)$ has isolated
vertices. Thus, the graph $G(\mathcal{X}_n;r_n)$ is connected w.h.p.

For $\alpha<\frac{2^{d-1}}{d\pi_d}$, by the proof of Theorem
\ref{result_0} (i) we see that for some vertex $k$,  w.h.p. there
exist some vertices lying in $B(x_k;n^{-1/d})$ which do not have
neighbors in $[0,1]^d\setminus B(x_k;n^{-1/d})$. Thus, the graph
$G(\mathcal{X}_n;r_n)$ is not connected w.h.p. This completes the
proof of the lemma.

\section*{Proof of Proposition \ref{Special}}\label{probosition}
 (i) By the definition of $\xi_n$, we have
\begin{eqnarray*}\label{s1}
\begin{aligned}
&E\left[{\rm e}^{\overline{\theta}_n \xi_n} \right]=1-\pi_d r_n^d
+d\pi_d
\int_0^{r_n}{\rm e}^{\overline{\theta}_n f_n(x)}x^{d-1}dx\\
&~~~~~~~~ ~~~~\,= 1-\pi_d r_n^d +d\pi_d r_n^d \int_0^{1}{\rm
e}^{\overline{\theta}_n
f_n(r_n y)}y^{d-1}dy,\\
&E\left[\xi_n {\rm e}^{\overline{\theta}_n \xi_n} \right]=d\pi_d
r_n^d \int_0^{1}f_n(r_n y) {\rm e}^{\overline{\theta}_n f_n(r_n
y)}y^{d-1}dy.
\end{aligned}
\end{eqnarray*}
Substituting the above two equations into (\ref{com4}) and using the
assumption  $\overline{\theta}_n=O(1/f_n(0))$, we have
\begin{eqnarray*}\label{s2}
&&n\overline{\theta}_n d\pi_d r_n^d \int_0^{1}f_n(r_n y)
{\rm e}^{\overline{\theta}_n f_n(r_n y)}y^{d-1}dy\\
&=&\left( \log n + n \left(-\pi_d r_n^d +d\pi_d r_n^d \int_0^{1}{\rm
e}^{\overline{\theta}_n f_n(r_n y)}y^{d-1}dy   \right)\right)(1\pm
O(r_n^d)),
\end{eqnarray*}
which indicates that
\begin{eqnarray}\label{s3}
\begin{aligned}
&\int_0^{1}\left(\overline{\theta}_nf_n(r_n y)-1\right) {\rm
e}^{\overline{\theta}_n f_n(r_n y)}y^{d-1}dy= \left(\frac{\log
n}{d\pi_d n r_n^d}-\frac{1}{d}\right)\left(1 \pm
O\left(r_n^d\right)\right) \pm O\left(r_n^d\right).
\end{aligned}
\end{eqnarray}
Let $\theta_n^*$ denote a solution of the following equation
\begin{eqnarray}\label{s3_1}
\begin{aligned}
\int_0^{1}\left(\theta_nf_n(r_n y)-1\right) {\rm e}^{\theta_n
f_n(r_n y)}y^{d-1}dy=\frac{\log n}{d\pi_d n r_n^d}-\frac{1}{d}
\end{aligned}
\end{eqnarray}
with respect to $\theta_n$. Note that if $\theta_n=0$ then
\bna \label{B1}\mbox{the left side of (\ref{s3_1}) $<$ the right side of (\ref{s3_1})},\ena
and if $\theta_n>\log n/(d\pi_d c_0 n r_n^d f_n(0))$, then by (\ref{S0_1}),
\bna\label{s4}
&&\mbox{the left side of (\ref{s3_1})}\nn\\
&=&\frac{1}{d}\left[\int_0^{1}\theta_nf_n\left(r_n z^{1/d}\right)
{\rm e}^{\theta_n f_n(r_n z^{1/d})}dz-\int_0^{1}{\rm e}^{\theta_n
f_n(r_n z^{1/d})}dz
\right]\nn\\
&\geq& \frac{1}{d}\left[\int_0^{1}\theta_nf_n\left(r_n
z^{1/d}\right)dz \int_0^1 {\rm e}^{\theta_n f_n(r_n
z^{1/d})}dz-\int_0^{1}{\rm e}^{\theta_n f_n(r_n z^{1/d})}dz
\right]\nn\\
&\geq& \frac{1}{d}\left(\theta_n
dc_0f_n(0)-1\right)\int_0^{1}{\rm e}^{\theta_n f_n(r_n z^{1/d})}dz\nn\\
&>&\mbox{the right side of (\ref{s3_1})},
\ena
where the following inequality is used
\begin{eqnarray*}
\begin{aligned}
&\frac{a_1b_1+\cdots+a_m b_m}{m}-\frac{a_1+\cdots a_m}{m}\cdot\frac{b_1+\cdots+b_m}{m}=\frac{1}{m^2}\sum_{1\leq i<j\leq m}(a_i-a_j)(b_i-b_j)\geq 0
\end{aligned}
\end{eqnarray*}
with $\{a_{i}, 1\leq i\leq m\}$ and $\{b_{i}, 1\leq i\leq m\}$ being two real number sequences satisfying $a_1\geq a_2\geq \cdots \geq a_m$ and $b_1\geq b_2\geq \cdots\geq b_m$.
Hence, by \dref{B1} and \dref{s4} we have
\begin{eqnarray}\label{s4_1}
0<\theta_n^*<\frac{\log n}{d\pi_d c_0 n r_n^d f_n(0)}.
\end{eqnarray}
By (\ref{s4_1}) and (\ref{rad_con}), we have $\theta_n^*=O(1/f_n(0))$.
Substituting this into (\ref{s3}), we see that the equation
(\ref{com4}) has a solution near to $\theta_n^*$. By the uniqueness
of the solution of (\ref{com4}) we obtain that $\overline{\theta}_n =
\theta_n^*(1+o(1))=O(1/f_n(0))$. Moreover, by the first equation of
(\ref{com4}), we have
\begin{eqnarray}\label{s5}
\overline{k}_n=n d\pi_d r_n^d \int_0^{1}f_n(r_n y) {\rm
e}^{\theta_n^* f_n(r_n y)}y^{d-1}dy(1+o(1))={\it \Theta}(nr_n^d
f_n(0)).
\end{eqnarray}

(ii) Similar to the analysis of (i), we have
\begin{eqnarray*}\label{s5_1}
\begin{aligned}
&E\left[{\rm e}^{\theta \xi_{n,\delta}} \right]= 1-\pi_d
r_n^d(1+\delta)^d +d\pi_d r_n^d (1+\delta)^d \int_0^{1}{\rm
e}^{\theta
f_{n,\delta}[(1+\delta)r_n y]}y^{d-1}dy,\\
&E\left[\xi_{n,\delta} {\rm e}^{\theta \xi_{n,\delta}}
\right]=d\pi_d r_n^d (1+\delta)^d
\int_0^{1}f_{n,\delta}[(1+\delta)r_n y] {\rm e}^{\theta
f_{n,\delta}[(1+\delta)r_n y]}y^{d-1}dy.
\end{aligned}
\end{eqnarray*}
Now, we consider the solution of the following equation,
\begin{eqnarray}\label{s6}
\int_0^{1}\left(\theta f_{n,\delta}[(1+\delta)r_n y]-1\right) {\rm
e}^{\theta f_{n,\delta}[(1+\delta)r_n y]}y^{d-1}dy=\frac{\log
n}{d\pi_d n r_n^d(1+\delta)^d}-\frac{1}{d}.
\end{eqnarray}
Define $g(y):=f_{n,\delta}[(1+\delta)r_n y]-f_n(r_n y)$. It is clear that $g(y)\geq 0$
for $y\leq 1$. Note that
\begin{eqnarray*}\label{s7}
\begin{aligned}
&\int_0^{1}\left(\theta_n^* f_{n,\delta}[(1+\delta)r_n y]-1\right)
{\rm e}^{\theta_n^* f_{n,\delta}[(1+\delta)r_n
y]}y^{d-1}dy-\frac{\log n}{d\pi_d n
r_n^d(1+\delta)^d}+\frac{1}{d}\\
>&\int_0^{1}\left\{\theta_n^* [f_{n}(r_n y)+g(y)]-1\right\}
{\rm e}^{\theta_n^* f_{n}(r_n y)} {\rm e}^{\theta_n^* g(y)}y^{d-1}dy\\
& -\int_0^{1}\left(\theta_n^* f_{n}(r_n y)-1\right)
{\rm e}^{\theta_n^* f_{n}(r_n y)}y^{d-1}dy\\
=&\int_0^{1}\left(\theta_n^* f_{n}(r_n y)-1\right) {\rm
e}^{\theta_n^*
f_{n}(r_n y)}\left({\rm e}^{\theta_n^* g(y)}-1 \right) y^{d-1}dy\\
& +\int_0^{1} \theta_n^* g(y){\rm e}^{\theta_n^* f_{n}(r_n y)}
{\rm e}^{\theta_n^* g(y)}y^{d-1}dy\\
>&\int_0^{1}{\rm e}^{\theta_n^* f_{n}(r_n y)}\left[\theta_n^*
g(y){\rm e}^{\theta_n^* g(y)}-{\rm e}^{\theta_n^* g(y)}+1 \right] y^{d-1}dy\\
=&\int_0^{1}{\rm e}^{\theta_n^* f_{n}(r_n
y)}\left[\sum_{i=2}^{\infty} \frac{i-1}{i!}(\theta_n^* g(y))^i
\right] y^{d-1}dy>0.
\end{aligned}
\end{eqnarray*}
If $\theta=0$, then
$$\mbox{the left side of (\ref{s6}) $<$ the right side of (\ref{s6})}.$$
Thus, the solution of equation (\ref{s6}) satisfies
$\theta_{n,\delta}^* \in (0,\theta_n^*)$. With a similar argument,
we have $\overline{\theta}_{n,\delta}=\theta_{n,\delta}^*(1+o(1))$.
For any $\delta>0$, similar to (\ref{s5})  we have
\ay\begin{eqnarray}\label{s8} \frac{\overline{k}_{n,\delta}}{n
d\pi_d r_n^d}&=& (1+\delta)^d\int_0^{1}f_{n,\delta}[(1+\delta)r_n y]
{\rm e}^{\theta_{n,\delta}^*
f_{n,\delta}[(1+\delta)r_n y]}y^{d-1}dy(1+o(1))\nonumber\\
&=&\int_0^{1+\delta}f_{n,\delta}(r_n z) {\rm e}^{\theta_{n,\delta}^*
f_{n,\delta}(r_n z)}z^{d-1}dz(1+o(1)).
\end{eqnarray}
By the definition of $f_{n,\delta}$ and the fact that
$\theta_{n,\delta}^*<\theta_n^*$ we have
\begin{eqnarray*}\label{s8_1}
\begin{aligned}
&\int_0^{1+\delta}f_{n,\delta}(r_n z) {\rm e}^{\theta_{n,\delta}^*
f_{n,\delta}(r_n z)}z^{d-1}dz\\
=&\int_0^{\delta}f_n(0) {\rm e}^{\theta_{n,\delta}^*
f_n(0)}z^{d-1}dz+ \int_0^{1}f_{n}(r_n z) {\rm
e}^{\theta_{n,\delta}^* f_{n}(r_n
z)}(z+\delta)^{d-1}dz\\
 <&\frac{\delta^d f_n(0) {\rm e}^{\theta_{n}^* f_n(0)} }{d}+
\sum_{i=0}^{d-1}\binom{d-1}{i}\delta^i \int_0^{1}f_{n}(r_n z)
{\rm e}^{\theta_{n}^* f_{n}(r_n z)}z^{d-1-i}dz\\
\leq& f_n(0) {\rm e}^{\theta_{n}^* f_n(0)}\left( \frac{\delta^d
}{d}+\sum_{i=1}^{d-1}\binom{d-1}{i}\frac{\delta^i
}{d-i}\right)+\int_0^{1}f_n(r_n z) {\rm e}^{\theta_n^* f_n(r_n
z)}z^{d-1}dz\\
=& (1+O(\delta))\int_0^{1}f_n(r_n z) {\rm e}^{\theta_n^* f_n(r_n
z)}z^{d-1}dz,
\end{aligned}
\end{eqnarray*}
where the condition (\ref{S0_1}) and the fact
$\theta_{n}^*=O(1/f_n(0))$ are used in the last equation. Combining
this with (\ref{s8}) and (\ref{s5}) our
result follows.\\
(iii) Immediate from (\ref{s4_1}) and (\ref{s5}).

\section*{Proof of Lemma \ref{Delta1}}\label{Delta1proof}

Before the proof of Lemma \ref{Delta1}, we need to introduce some
notations. For any $n\in \mathds{N}$, define
$\rho(n):=\lfloor n-2n^{3/4}\rfloor$ and
\begin{eqnarray*}\label{a1_4}
\widehat{I}_{n,\delta}(x):=\sup_{\theta>0}\left\{\theta
x-(\rho(n)-1) \log \left(E\left[{\rm e}^{\theta \xi_{n,\delta} }
\right]\right) \right\}.
\end{eqnarray*}
Take $\widehat{k}_{n,\delta}>(\rho(n)-1)
E[\xi_{n,\delta}]$ such that
$\widehat{I}_{n,\delta}(\widehat{k}_{n,\delta})=\log n$. Denote
$\widetilde{f}_n=\widetilde{f}_{n,0}$,
$\widetilde{\xi}_{n}=\widetilde{\xi}_{n,0}$,
$\widetilde{I}_{n}=\widetilde{I}_{n,0}$ and
$\widehat{k}_{n}=\widehat{k}_{n,0}$.

\begin{lemma}\label{lem_2} For $\widehat{k}_{n}$ and $\overline{k}_n$, we have
 $
\lim_{n\rightarrow\infty}( \widehat{k}_{n}/\overline{k}_n )=1.$
\end{lemma}
\proof First, by  $$ \widehat{I}_n(\overline{k}_n) \geq
I_n(\overline{k}_n)=\log n=\widehat{I}_n(\widehat{k}_n),$$ the
inequality $\overline{k}_n\geq \widehat{k}_n $ can be derived.

Assume that there exists a constant $\varepsilon>0$ such that
$\overline{k}_n\geq(1+\varepsilon)\widehat{k}_n$. Then for large
$n$,
\begin{eqnarray*}
\begin{aligned}
I_n(\overline{k}_n)&\geq I_n\left((1+\varepsilon)\widehat{k}_n
\right)\\
&=\sup_{\theta>0}\left\{\theta (1+\varepsilon)\widehat{k}_n-(n-1)
\log \left(E\left[{\rm e}^{\theta\xi_n } \right]\right) \right\}\\
&\geq\frac{n-1}{\rho(n)-1}\sup_{\theta>0}\left\{\theta
\left(1+\frac{\varepsilon}{2}\right)\widehat{k}_n-(\rho(n)-1) \log
\left(E\left[{\rm e}^{\theta\xi_n } \right]\right) \right\}\\
&\geq\frac{n-1}{\rho(n)-1} \left(1+\frac{\varepsilon}{2}\right)
\widehat{I}_n(\widehat{k}_n)\\
&>\log n=I_n(\overline{k}_n),
\end{aligned}
\end{eqnarray*}
which leads to contradiction. Our result yields.  \rulex

\emph{Proof of  Lemma} \ref{Delta1}~ For simplicity of expressions,
we consider the case of $\delta=0$, and it is easy to extend our
results to the case of $\delta>0$.

Denote $d_i=d_i(n):=\sum_{1\leq j \leq n, j\neq
i}f_n(\|X_i(0)-X_j(0)\|)$. For $k_n>(n-1) E[\xi_n]$, by (2.2.12) in
\cite{Dem1998}, we have \ay\begin{eqnarray}\label{d1}
P(d_i>k_n)&\leq& P\left(\sum_{1\leq j \leq n, j\neq
i}f_n\left(\big\|x_0-X_j(0)\big\|\right)>k_n
\right)\nonumber\\
&=&P\left(\frac{1}{n-1}\sum_{1\leq j \leq n, j\neq
i}f_n\left(\big\|x_0-X_j(0)\big\|\right)>\frac{k_n}{n-1}
\right)\nonumber\\
&\leq &\exp\left(-(n-1)\sup_{\theta>0}\left\{\frac{\theta k_n}{n-1}
-\log \left(E\left[{\rm e}^{\theta\xi_{n} } \right]\right)
\right\}\right)={\rm e}^{-I_n(k_n)},
\end{eqnarray}
where $I_n(\cdot)$ is defined by (\ref{com1}). For any
$\varepsilon>0$ we have
\bna\label{d5}
&&I_n((1+\varepsilon)\overline{k}_n)\nn\\&\geq& \overline{\theta}_n
\overline{k}_n(1+\varepsilon)-(n-1)
\log \left(E\left[{\rm e}^{\overline{\theta}_n\xi_n } \right]\right) \ \ \ \ \hbox{(according to (\ref{com1})} )\nn\\
&\geq& (1+\varepsilon)\left( \overline{\theta}_n
\overline{k}_n-(n-1) \log \left(E\left[{\rm
e}^{\overline{\theta}_n\xi_n
} \right]\right)  \right)\ \ \ \ \hbox{(according to (\ref{com4}))}\nn\\
&=& (1+\varepsilon)\log n.
\ena
Combining \dref{d5} with (\ref{d1}), we obtain
\begin{eqnarray}\label{d6}
P(d_i>(1+\varepsilon)\overline{k}_n)\leq {\rm
e}^{-I_n((1+\varepsilon)\overline{k}_n)}\leq {\rm
e}^{-(1+\varepsilon)\log n}=n^{-1-\varepsilon}.
\end{eqnarray}

Set $$\overline{F}_n:=\bigcup_{i=1}^n
I_{\{d_i>(1+\varepsilon)\overline{k}_n\}}.$$ By the Boole's
inequality and (\ref{d6}), we have
\begin{eqnarray*}\label{d7}
P(\overline{F}_n)\leq \sum_{i=1}^n
P(d_i>(1+\varepsilon)\overline{k}_n) \leq n^{-\varepsilon}.
\end{eqnarray*}
Thus,
\begin{eqnarray*}
P\left(\Delta_n>(1+\varepsilon)\overline{k}_n\right)=
P(\overline{F}_n)\rightarrow 0~~~~\mbox{as}~~n\rightarrow\infty.
\end{eqnarray*}
The inequality (\ref{D1_00}) holds.  \rulex

\section*{Proof of Lemma \ref{Eigenvalue}}\label{eigenvalueproof}

The proof of Lemma \ref{Eigenvalue} mainly uses the idea appearing
in the proof of Theorem 4.3 of \cite{Chen2014}.
We first introduce some notations.

For any constant $\varepsilon\in (0,1)$, set
$K_n=K_n(\varepsilon):=\lceil\frac{\sqrt{d+3}}{\varepsilon
r_n}\rceil$, where $\lceil x\rceil$ is the smallest integer no less
than $x$. The unit square $[0,1]^d$ is divided into $K_n^d$ equal
small squares with the length of each side equal to $1/K_n$. We
denote these small squares as $S_1(n),S_2(n),\cdots,S_{K_n^d}(n)$.
Denote $S_i=S_i(n)$ for $1\leq i \leq K_n^d$. For each small square
$S_i$, $1\leq i \leq K_n^d$, let $x_i$ denote its center point, and
$z_i:=K_n x_i+(\frac{1}{2},\frac{1}{2},\cdots,\frac{1}{2}) \in
\mathbb{Z}^d$.

Let $\|\cdot\|_1$ and $\|\cdot\|_{\infty}$ denote the $l_1$-norm and
infinity norm respectively. For any $x,y\in \mathbb{Z}^d$, if
$\|x-y\|_1=1$, then we say that $x$ and $y$ are  {adjacent}, written
as $x\sim y$. Given $A\subseteq \mathbb{Z}^d$, if for any $x,y\in
A$, there exists a vertex sequence $x_1,x_2,\cdots,x_n$ in $A$ such
that $x\sim x_1, x_1\sim x_2, x_2\sim x_3,\cdots, x_n\sim y$, then
we say $A$ is connected. Similarly, if $\|x-y\|_{\infty}\leq k$,
$k\geq 1$, we say that $x$ and $y$ are {$k$-adjacent}, written as
$x\sim_k y$. Given $A\subseteq \mathbb{Z}^d$, if for any $x,y\in A$,
there exists a vertex sequence $x_1,x_2,\cdots,x_n$ in $A$ such that
$x\sim_k x_1, x_1\sim_k x_2, x_2\sim_k x_3,\cdots, x_n\sim_k y$,
then we say $A$ is $k$-{connected}. We see that for any $k\geq 1$ if
$A$ is $k$-connected, then $A$ must be connected. In particular, a
single vertex set $\{x\}\subset \mathbb{Z}^d$ is both connected and
$k$-connected.

We define the lattice box $B_{\mathbb{Z}}(K_n)$ by
$B_{\mathbb{Z}}(K_n):=\prod_{i=1}^d([1,K_n]\cap \mathbb{Z}).$ It is
clear that $B_{\mathbb{Z}}(K_n)$ is equal to the set $\{z_i:1\leq i
\leq K_n^d\}$. For $A\subset B_{\mathbb{Z}}(K_n)$, we denote
$A^c:=B_{\mathbb{Z}}(K_n)\backslash A$. Let $\partial A$ denote the
 {internal vertex-boundary} of $A$, that is, the set of vertex $z\in
A$ such that $\{y\in A^c : \|z-y\|_1=1\}$ is non-empty.

For $\eta>0$, we use $Po(\eta)$ to denote the Poisson random
variable with parameter $\eta$. Define a Poisson point process
$\mathcal{P}_{\eta}$ as
$\mathcal{P}_{\eta}:=\{Y_1,Y_2,\cdots,Y_{Po(\eta)}\}$, where
$\{Y_1,Y_2,\cdots\}$ is the set of vertices independently and
uniformly distributed in $[0,1]^d$ and $Po(\eta)$ is independent of
$\{Y_1,Y_2,\cdots\}$, see Subsection 1.7 in \cite{penrose2003}. For
a Borel set $A\subseteq [0,1]^d$, $|\mathcal{P}_{\eta}\cap A|$ is a
Poisson random variable with parameter $\eta Leb(A)$, where
$|\cdot|$ denotes the cardinality and $Leb(\cdot)$ denotes the
Lebesgue measure. For any two Borel sets $A_1,A_2\subseteq [0,1]^d$,
if $Leb(A_1 \cap A_2)=0$, then the random variables
$|\mathcal{P}_{\eta} \cap A_1|$ and $|\mathcal{P}_{\eta} \cap A_2|$
are mutually independent. Set $\eta(n):=n-n^{3/4}$, and let
$\mathcal{P}_{\eta(n)}$ be a Poisson point process in $[0,1]^d$ with
parameter $\eta(n)$. Then, $\mathcal{P}_{\eta(n)}\subseteq
\mathcal{X}_n$ except when $Po(\eta(n))>n$, and by Lemma 1.4 in
\cite{penrose2003} we obtain
\begin{eqnarray}\label{d13}
 \mathcal{P}_{\eta(n)}\subseteq
\mathcal{X}_n~~~~~~~\mbox{w.h.p.}
\end{eqnarray}

For any set $A\subseteq B_{\mathbb{Z}}(K_n)$,  define the function
\begin{eqnarray*}
g_1(A):=\sum_{z_i\in A,z_j\in A^c,z_i\sim z_j} |
\mathcal{P}_{\eta(n)}\cap S_i|\cdot | \mathcal{P}_{\eta(n)}\cap
S_j|.
\end{eqnarray*}

For any $z_i \in B_{\mathbb{Z}}(K_n)$, we call $z_i$  {open} if $S_i
\cap \mathcal{P}_{\eta(n)} \neq \emptyset$, and call $z_i$
\emph{closed} otherwise. Let $\mathcal{O}_n$ denote the set of open
vertices in $B_{\mathbb{Z}}(K_n)$, and let $\mathcal{C}_n$ denote
the largest open clusters of $\mathcal{O}_n$.

Before proving Lemma \ref{Eigenvalue}, we provide some preliminary
results, see the following Lemmas{}\,\ref{Delta}--\ref{Cheeger_lem}.
 Lemmas \ref{Delta} and \ref{iso2}--\ref{Cheeger_lem} are proved
under the condition (\ref{rad_con}), and we will omit it to avoid
repetition.

\begin{lemma}\label{Delta}
There exists a constant $c=c(\varepsilon,\alpha,d)$ such that
\begin{eqnarray*}
\max_{1\leq i \leq K_n^d}|\mathcal{X}_n\cap S_i|\leq c \varepsilon^d
n r_n^d~~~~~~~~\mbox{w.h.p.}
\end{eqnarray*}
\end{lemma}
\proof
It can be  easily deduced from Lemma 4.1 of \cite{Chen2014}.  \rulex

\begin{lemma}[Lemma 9.9 in \cite{penrose2003}]\label{Isoperimetric}
Let $\beta \in (0,1)$. If $A$ is a subset of $B_{\mathbb{Z}}(K_n)$
$($not necessarily connected$)$, with $|A|\leq \beta K_n^d$, then
\begin{eqnarray*}\label{I0}
|\partial A|\geq (2d)^{-1}(1-\beta^{1/d}){|A|}^{(d-1)/d}.
\end{eqnarray*}
\end{lemma}

\begin{lemma}\label{iso2}
Suppose that $A\subset B_{\mathbb{Z}}(K_n)$ and
the integer $k\geq 1$. Then for any $\beta\in(0,1)$, there
exist constants $c=c(\alpha,\varepsilon,k,\beta,d)>0$ and $\gamma=\gamma(\alpha,\varepsilon,k,\beta)$ such that w.h.p.
\begin{eqnarray*}
\inf_{\gamma (\log n)^{d/(d-1)} \leq |A| \leq \beta K_n^d \atop A
\mbox{\rm \small is $k$-connected }} \frac{g_1(A)}{|A|}>  \frac{cn^2
r_n^{2d}}{K_n}.
\end{eqnarray*}
\end{lemma}
\proof
The lemma can be deduced by the similar method as that of Lemma 5.10
in \cite{Chen2014} with a small modification, and we omit the proof details to save space.  \rulex

\begin{lemma}\label{limitinfty}
If $\alpha \varepsilon^d > (d+3)^{d/2}$ then w.h.p.
$$\min_{1\leq i\leq K_n^d}|S_i \cap \mathcal{P}_{\eta(n)}| > \frac{n}{K_n^d} H_-^{-1}\left(\frac{1}{2}+\frac{(d+3)^{d/2}}{2\alpha \varepsilon^d} \right).$$
\end{lemma}
\proof
By Lemma 1.2 in
\cite{penrose2003} we obtain
\begin{eqnarray*}
\begin{aligned}
&P\left(\bigcup_{i=1}^{K_n^d} \left\{|\mathcal{P}_{\eta(n)}\cap
S_{i}|\leq \frac{n\beta}{K_n^d}\right\}\right)\\  \leq&
K_n^d\exp\left(-\frac{n-n^{3/4}}{\lceil
\frac{\sqrt{d+3}}{\varepsilon r_n} \rceil^d } H\left( \frac{n\beta}{n-n^{3/4}} \right) \right)\\
<& n \exp\left(-\frac{n-n^{3/4}}{\lceil
\frac{\sqrt{d+3}}{\varepsilon r_n} \rceil^d }
\left(\frac{1}{3}+\frac{2}{3 \alpha\varepsilon^d(d+3)^{-d/2}}
\right) \right)\rightarrow 0~~\mbox{as}~~ n\rightarrow\infty,
\end{aligned}
\end{eqnarray*}
which is followed by our result.  \rulex

\begin{lemma}\label{iso3}
Suppose that
$\alpha \varepsilon^d > (d+3)^{d/2}$  and $A\subset B_{\mathbb{Z}}(K_n)$. Then
for any $\beta\in(0,1)$, there exists a constant
$c=c(\alpha,d,\beta)>0$ such
that w.h.p.
\begin{eqnarray*}
\inf_{ |A| \leq \beta K_n^d} \frac{g_1(A)}{|A|}>c
n^2 r_n^{2d+1}.
\end{eqnarray*}
\end{lemma}
\proof
Since $|A| \leq \beta K_n^d$, then by Lemma \ref{Isoperimetric}, we have
\begin{eqnarray}\label{I1}
|\partial A|\geq (2d)^{-1}(1-\beta^{1/d}){|A|}^{(d-1)/d} \geq
\frac{(1-\beta^{1/d})|A|}{2d\beta^{1/d} K_n}.
\end{eqnarray}
Combining \dref{I1} with Lemma \ref{limitinfty}, w.h.p.
\begin{eqnarray*}
\begin{aligned}
\inf_{|A| \leq \beta K_n^d} \frac{g_1(A)}{|A|} &>
\frac{n^2}{K_n^{2d}}
\left(H_-^{-1}\left(\frac{1}{2}+\frac{(d+3)^{d/2}}{2\alpha \varepsilon^d} \right)\right)^2
\inf_{|A| \leq \beta K_n^2}
\frac{|\partial A|}{|A|}\\
& \geq \frac{ n^2 (1-\beta^{1/d})}{\beta^{1/d} K_n^{2d+1}}
\left(H_-^{-1}\left(\frac{1}{2}+\frac{(d+3)^{d/2}}{2\alpha\varepsilon^d} \right)\right)^2,
\end{aligned}
\end{eqnarray*}
which implies our result.  \rulex

~~ 

For $F\subseteq \{1,2,\cdots,n\}$ and $F^c=\{1,2,\cdots,n\}
\backslash F$, set
$$\Phi_n:=\inf_{|F|\leq \frac{n}{2}, t>0}\left\{\frac{1}{|F|}\sum_{i\in F, j\in F^c}
f_n(\|X_i(t)-X_j(t)\|)\right\}.$$
We have the following lemma.

\begin{lemma}\label{Cheeger_lem}
Assume that $\alpha \in (\frac{2^{d-1}}{d\pi_d},\infty]$, and that
there exists a constant $\delta>0$ such that{}\,$(\ref{Ei_0})$
holds. For any constant $\varepsilon
>0$,  there exists a
constant $c=c(\alpha,\varepsilon,d)>0$ such that

{\rm (i)} if $\alpha\varepsilon^d\leq (d+3)^{d/2} $, then w.h.p.
\begin{eqnarray*}
{\it \Phi}_n \geq c \min \left\{nr_n^{d+1}
f_n((\delta+\varepsilon)r_n), \frac{
f_n\left(R_c+(\delta+\varepsilon)r_n \right)}{(\log
n)^{(2d-1)/(d-1)}} \right\};
\end{eqnarray*}

{\rm (ii)} if $\alpha\varepsilon^d > (d+3)^{d/2}$, then w.h.p.
\begin{eqnarray*}
{\it \Phi}_n \geq  c n f_n((\delta+\varepsilon)r_n)\min\left\{
r_n^{d+1},1\right\}.
\end{eqnarray*}
\end{lemma}

\proof (i)~ For  $F\subseteq \{1,2,\cdots,n\}$, define
$\widetilde{F}:=\{X_i(0):i\in F\}\subseteq \mathcal{X}_n$ to be the
initial positions of agents whose indexes are in $F$. For
$D_1,D_2\subset [0,1]^d$, set
\begin{eqnarray*}
g_{D_1,D_2}(F)=g_{D_1,D_2,n}(F):=\sum_{x\in D_1\cap\widetilde{F},
y\in D_2\cap \widetilde{F^c}}f_n(\|x-y\|+\delta r_n).
\end{eqnarray*}
Take $g(F)=g_{[0,1]^d,[0,1]^d}(F)$, and define
\begin{eqnarray*}
\begin{aligned}
{\it \Phi}_n'&:=\inf_{|F|\leq n/2}\frac{\sum_{i\in F, j\in F^c}
f_n(\|X_i(0)-X_j(0)\|+\delta
 r_n)}{|F|}=\inf_{|F|\leq n/2}
\frac{g(F)}{|F|}.
\end{aligned}
\end{eqnarray*}

Now, we estimate ${\it \Phi}_n'$. Denote
$$A_F:=\left\{z_i: |S_i \cap \widetilde{F}|> \frac{1}{2} |S_i\cap \mathcal{X}_n|\right\}\subseteq
B_{\mathbb{Z}}(K_n) $$ and
\begin{eqnarray*}
\widetilde{A_F}:=\bigcup_{z_i\in A_F} S_i\cap \mathcal{X}_n.
\end{eqnarray*}
Set $\beta:=1-\frac{1}{4c_1(d+4)^{d/2}}$. If $|A_F|>\beta K_n^d$,
then $|A_F^c|\leq (1-\beta) K_n^d$. By Lemma \ref{Delta}, we have w.h.p.
$$\sum_{z_i\in A_F^c} |S_i \cap \mathcal{X}_n| \leq c_1\varepsilon^d|A_F^c|
n r_n^d < \frac{n}{4},$$ where $c_1$ is the constant $c$  appearing in Lemma \ref{Delta}. If $|F|\leq n/2$,
then $|\widetilde{F^c}|=|F^c|>\frac{n}{2}$. Thus, there exist at least
$\frac{n}{4}$ vertices in $\widetilde{F^c}$ contained by
$\widetilde{A_F}$. For $x\in \widetilde{F^c} \cap \widetilde{A_F}$,
without loss of generality we assume that $x\in S_i$ with $z_i\in
A_F$. Then by the definition of $A_F$ we can get $|\widetilde{F}\cap
S_i|\geq |\widetilde{F^c}\cap S_i| \geq 1$, which indicates that
there exists at least one vertex $y$ such that $y\in \widetilde{F}\cap S_i$. Note
that the length of the side of $S_i$ is less than $\varepsilon
r_n/\sqrt{d+3}$. Then w.h.p.
\begin{eqnarray}\label{E2}
\begin{aligned}
\inf_{|F| \leq \frac{n}{2},|A_F|>\beta K_n^d }\frac{g(F)}{|F|}\geq
\frac{\frac{n}{4}f_n\left(\frac{\varepsilon r_n
\sqrt{d}}{\sqrt{d+3}}+\delta r_n\right)}{n/2}\geq
\frac{f_n((\varepsilon+\delta)r_n)}{2}.
\end{aligned}
\end{eqnarray}

In the following we consider the case of $|A_F|\leq \beta K_n^d$.
Let $A_1,A_2,\cdots,A_{m_F}$ be components of $A_F$ satisfying:  1)
$A_1,A_2,\cdots,A_{m_F}$ are all
$\lceil\frac{\sqrt{d+3}}{\varepsilon}\rceil$-connected;  2) $A_i\cup
A_j$, $1\leq i\neq j\leq m_F$ is not
$\lceil\frac{\sqrt{d+3}}{\varepsilon}\rceil$-connected;  3)
$|A_1|\geq |A_2|\geq\cdots\geq A_{m_F}$. Without loss of generality,
we assume that $|A_i|\geq \gamma(\log n)^{d/(d-1)}$ for $1\leq i
\leq i_F$, and $|A_i| < \gamma(\log n)^{d/(d-1)}$ for $i_F+1\leq i
\leq m_F$, where $i_F\in [1,m_F]$ and $\gamma$ is the same constant
appearing in Lemma \ref{iso2}. By Lemma \ref{iso2}, we have
\begin{eqnarray}\label{E4}
\inf_{|A_F|\leq \beta K_n^d,1\leq i \leq i_F}
\frac{g_1(A_i)}{|A_i|}\geq \frac{c_2 n^2 r_n^{2d}}{K_n}~~~~\mbox{w.h.p.},
\end{eqnarray}
where $c_2$ is the constant $c$ appearing in Lemma
\ref{iso2}.

For $i\in [1,i_F]$,  it is easy to see that if
 $z_k\in A_i$ and $z_j\in A_i^c$ with $z_k\sim z_j$, then $z_j\in A_F^c$, and the distance of any pair of vertices in
$S_k\cup S_j$ is not greater than $\varepsilon r_n$. By the
definition of $A_F$, we have
\begin{eqnarray*}\label{E5}
\begin{aligned}
g_{S_k,S_j}(F)&=\sum_{x\in S_k\cap \widetilde{F},y\in S_j\cap
\widetilde{F^c}}f_n(\|x-y\|)\geq \frac{f_n((\varepsilon+\delta) r_n)}{4}|\mathcal{X}_n\cap
S_k|\cdot |\mathcal{X}_n\cap S_j|.
\end{aligned}
\end{eqnarray*}
Therefore, if $\mathcal{P}_{\theta(n)}\subseteq \mathcal{X}_n$, then
\begin{eqnarray}\label{E6}
\begin{aligned}
\sum_{z_k\in A_i,z_j\in A_F^c,z_k\sim
z_j}g_{S_k,S_j}(F)&=\sum_{z_k\in A_i,z_j\in A_i^c,z_k\sim
z_j}g_{S_k,S_j}(F)\geq \frac{f_n((\varepsilon+\delta) r_n)}{4}g_1(A_i).
\end{aligned}
\end{eqnarray}
Moreover, by (\ref{d13}), we see that
$\mathcal{P}_{\theta(n)}\subseteq \mathcal{X}_n$ holds w.h.p.

Set
\begin{eqnarray*}
S_F^1:=\bigcup_{i=1}^{i_F} \bigcup_{z_k\in A_i} S_k.
\end{eqnarray*}
By (\ref{E6}), we have w.h.p.
\begin{eqnarray}\label{E7}
\begin{aligned}
 g_{S_F^1,[0,1]^d\backslash
S_F^1}(F)&\geq\sum_{i=1}^{i_F} \sum_{z_k\in A_i,z_j\in A_F^c,z_k\sim
z_j}g_{S_k,S_j}(F)\geq \sum_{i=1}^{i_F}\frac{f_n((\varepsilon+\delta)
r_n)}{4}g_1(A_i).
\end{aligned}
\end{eqnarray}

For $i\in [i_F+1,m_F]$, if $\bigcup_{z_j\in A_i} S_j\cap
\widetilde{F^c}\neq \emptyset$, then we have $g_{D_i',D_i'}(F)\geq f_n(\varepsilon
r_n)$ where $D_i'=\bigcup_{z_j\in
A_i} S_j$; If $\bigcup_{z_j\in A_i} S_j\cap \widetilde{F^c}=\emptyset$,
then by Lemma \ref{rgg}, we know that w.h.p.
$G(\mathcal{X}_n;R_c+\varepsilon r_n)$ is connected.
Thus, there exists at least one vertex $x^*\in (\bigcup_{z_j\in A_i}
S_j)^c\cap \mathcal{X}_n$ such that the set
$$\{y:y\in \bigcup_{z_j\in A_i} S_j\cap \widetilde{F}, \|x^*-y\|\leq
R_c+\varepsilon r_n\}$$ is not empty. Assume that $x^*\in S_k\
(1\leq k \leq K_n^d)$ and $z_k$ is the corresponding integer point
of $S_k$, then $z_k$ must be
$\lceil\frac{\sqrt{d+3}}{\varepsilon}\rceil$-connected with $A_i$,
and $z_k \in A_F^c$. Denote $D_i''=\bigcup_{z_j\in A_i} S_j \cup
S_k$. If $x^*\in \widetilde{F^c}$, then
\begin{eqnarray*}
g_{D_i'',D_i''}(F)\geq f_n(R_c+(\delta+\varepsilon)r_n);
\end{eqnarray*}
Otherwise, by the definition of $A_F$ we have $S_k\cap
\widetilde{F^c}\neq \emptyset$, and
$$g_{D_i'',D_i''}(F)\geq g_{S_k,S_k}(F)\geq f_n((\varepsilon+\delta)r_n).$$

Set
\begin{eqnarray*}
S_F^2:=\left\{%
\begin{array}{lll}
\d \bigcup_{i=i_F+1}^{m_F} D_i',~~~~~~~~&\mbox{if  $\bigcup_{z_j\in
A_i} S_j\cap
\d \widetilde{F^c}\neq \emptyset$},\\
\d \bigcup_{i=i_F+1}^{m_F} D_i'',~~~~~~~&\mbox{otherwise}.
\end{array}%
\right.
\end{eqnarray*}
For $z\in \mathbb{Z}^d$,
it is easy to see that the number of different
$\lceil\frac{\sqrt{d+3}}{\varepsilon}\rceil$-connected components
which is $\lceil\frac{\sqrt{d+3}}{\varepsilon}\rceil$-connected with
$z$ is less than $\lceil\frac{\sqrt{d+3}}{\varepsilon}\rceil^d$. By
the above argument we have w.h.p.
\begin{eqnarray}\label{E8}
g_{S_F^2,S_F^2}(F)\geq
\bigg\lceil\frac{\sqrt{d+3}}{\varepsilon}\bigg\rceil^{-d}(m_F-i_F)f_n\left(R_c+(\delta+\varepsilon)r_n\right).
\end{eqnarray}

Let $S_F^3=[0,1]^d\backslash (S_F^1\cup S_F^2)$. For $x\in S_F^3\cap
\widetilde{F}$, we assume that $x\in S_k\ (1\leq k\leq K_n^d)$, and
denote $z_k\in B_{\mathbb{Z}}(K_n)$ as the corresponding integer
point of $S_k$. It is clear that $z_k\in A_F^c$, and  the set
$S_k\cap \widetilde{F^c}$ is not empty. Thus,
\begin{eqnarray}\label{E9}
g_{S_F^3,S_F^3}(F)\geq \sum_{x\in S_F^3\cap
\widetilde{F}}f_n((\delta+\varepsilon)r_n) = |S_F^3\cap
\widetilde{F}|f_n((\delta+\varepsilon)r_n).
\end{eqnarray}

By the definition of $S_F^1$ and $S_F^2$ we have $Leb(S_F^1\cap
S_F^2)=0$. By (\ref{E7}), (\ref{E8}) and (\ref{E9}), we have w.h.p.
\begin{eqnarray*}\label{E10}
\begin{aligned}
g(F)&\geq  g_{S_F^1,[0,1]^d\backslash S_F^1}(F)+ g_{S_F^2,S_F^2}(F)+
g_{S_F^3,S_F^3}(F)\\
&\geq \sum_{i=1}^{i_F}\frac{f_n((\varepsilon+\delta)r_n)}{4}g_1(A_i)
+|S_F^3\cap \widetilde{F}|f_n((\delta+\varepsilon)r_n)\\
&~~~~+\bigg\lceil\frac{\sqrt{d+3}}{\varepsilon}\bigg\rceil^{-d}(m_F-i_F)f_n\left(R_c+(\delta+\varepsilon)r_n\right).
\end{aligned}
\end{eqnarray*}
By the above inequality, we have w.h.p.
\begin{eqnarray*}\label{E11}
\begin{aligned}
\inf_{|A_F|\leq \beta K_n^d }\frac{g(F)}{|F|}&= \inf_{|A_F|\leq
\beta K_n^d }\frac{g(F)}{|S_F^1\cap \widetilde{F}|+|S_F^2\cap
\widetilde{F}|+|S_F^3\cap \widetilde{F}|}\\
&\geq \inf_{|A_F|\leq \beta K_n^d } \frac{g(F)}{c_1\varepsilon^d n
r_n^d \left(\sum_{i=1}^{i_F}|A_i| +(m_F-i_F)\gamma(\log
n)^{d/(d-1)}\right) +|S_F^3\cap
\widetilde{F}|}\\
&\geq
\min\left\{\frac{\frac{f_n((\varepsilon+\delta)r_n)}{4}\sum_{i=1}^{i_F}g_1(A_i)}{c_1
\varepsilon^d n r_n^d \sum_{i=1}^{i_F}|A_i| },
\frac{\lceil\frac{\sqrt{d+3}}{\varepsilon}\rceil^{-d}f_n\left(R_c+(\delta+\varepsilon)r_n\right)}{c_1
\varepsilon^d n r_n^d \gamma(\log n)^{d/(d-1)}},f_n((\delta+\varepsilon)r_n)   \right\}\\
&\geq \min\left\{\frac{f_n((\varepsilon+\delta)r_n)c_2 n^2 r_n^{2d}}{4c_1 \varepsilon^d n r_n^d K_n},
\frac{f_n\left(R_c+(\delta+\varepsilon)r_n\right)}{\lceil\frac{\sqrt{d+3}}{\varepsilon}\rceil^{d}c_1
\varepsilon^d n r_n^d \gamma(\log n)^{d/(d-1)}} \right\},
\end{aligned}
\end{eqnarray*}
where (\ref{E4}) is used in the last inequality. Combining this with
(\ref{E2}), we obtain that there exists a constant $c>0$ such that
\begin{eqnarray}\label{E12}
{\it \Phi}_n' \geq c \min \left\{nr_n^{d+1}
f_n((\delta+\varepsilon)r_n), \frac{
f_n\left(R_c+(\delta+\varepsilon)r_n \right)}{(\log
n)^{(2d-1)/(d-1)}} \right\}~~~~\mbox{w.h.p.}
\end{eqnarray}

By (\ref{Ei_0}), we have for $t>0$,
\begin{eqnarray*}\label{E13}
\|X_i(t)-X_j(t)\| \leq \|X_i(0)-X_j(0)\| +\delta r_n,\quad
\mbox{w.h.p.}
\end{eqnarray*}
Hence, ${\it \Phi}_n \geq {\it \Phi}_n'$ holds w.h.p. Combining this
with (\ref{E12}) yields our result.

(ii)~The definition of $c_1$, $\beta$ and  $A_F$ is the same as that
in (i). By (\ref{E2}) and Lemma \ref{limitinfty}, we have
\ay\begin{eqnarray}\label{E14} \inf_{|F| \leq
\frac{n}{2},|A_F|>\beta K_n^d }\frac{g(F)}{|F|}&\geq&
\frac{\frac{n}{4}\cdot\frac{n}{2 K_n^d}
H_-^{-1}\left(\frac{1}{2}+\frac{(d+3)^{d/2}}{2\alpha\varepsilon^d}
\right)f_n\left(\frac{\varepsilon r_n \sqrt{d}}{\sqrt{d+3}}+\delta
r_n\right)}{n/2}\nonumber\\
&\geq &\frac{n
H_-^{-1}\left(\frac{1}{2}+\frac{(d+3)^{d/2}}{2\alpha\varepsilon^d}
\right) f_n((\varepsilon+\delta)r_n)}{4 K_n^d}~~~~\mbox{w.h.p.}
\end{eqnarray}
For the case of $|A_F|\leq \beta K_n^d$, by Lemmas
\ref{Delta} and \ref{iso3} we have
\begin{eqnarray*}
\begin{aligned}
\inf_{|F| \leq \frac{n}{2},|A_F|\leq \beta K_n^d }\frac{g(F)}{|F|}
&\geq \frac{f_n((\varepsilon+\delta)r_n)}{c_1\varepsilon^d n
r_n^d}\inf_{|A_F|\leq
\beta K_n^d }\frac{g_1(F)}{|A_F|}\geq \frac{c_3 n^2 r_n^{2d+1}}{c_1 \varepsilon^d n r_n^d }
f_n((\varepsilon+\delta)r_n)~~\mbox{w.h.p.},
\end{aligned}
\end{eqnarray*}
where $c_3$ is the same constant as $c$ appearing in Lemma
\ref{iso3}. By (\ref{E14}) and the fact ${\it \Phi}_n \geq {\it
\Phi}_n'$ yields our result. \rulex

\vskip 3mm

\emph{Proof of Lemma} \ref{Eigenvalue}~ Suppose
$\overline{k}_{n,\delta} \leq 1-\varepsilon$. By (\ref{Ei_0}), we
have for $t>0$,
\begin{eqnarray*}
\|X_i(t)-X_j(t)\| \geq \|X_i(0)-X_j(0)\|-\delta r_n, \quad
\mbox{w.h.p.}
\end{eqnarray*}
By (\ref{newsec_3}) $\Delta_{n}(t)\leq \Delta_{n,\delta}$
holds w.h.p. for all $t>0$. By Lemma \ref{Delta1}, we see
that w.h.p. $\Delta_{n}(t) \leq 1-\frac{2\varepsilon}{3}$ holds for all $t>0$. Thus, given $\lambda\in
\mathbb{R}$, if $\lambda<\frac{\varepsilon}{2}-1$, then w.h.p. $P(t)-\lambda I_n$ is a strictly
diagonally dominant matrix and $\det(P(t)-\theta I_n)\neq 0$ for all
$t>0$, which indicates that $\lambda$ is not an
eigenvalue of $P(t)$. Thus, w.h.p.
\begin{eqnarray}\label{Ei2_1}
\lambda_n(t)\geq\frac{\varepsilon}{2}-1,~~~~~~~~\forall~t>0.
\end{eqnarray}

On the other hand, note that $P(t)$ is a symmetric stochastic
matrix, then the stationary distribution of $P(t)$ is
$(\frac{1}{n},\frac{1}{n},\cdots,\frac{1}{n})$. Therefore for $t>0$,
the Cheeger' constant of $P(t)$ is not less than ${\it \Phi}_n$. By
the Cheeger's inequality (Proposition 6 in \cite{Diaconis1991}), we
have $\lambda_2(t) \leq 1-  {\it \Phi}_n^2$ for $t>0$. Combining
this with Lemma \ref{Cheeger_lem} and (\ref{Ei2_1}), our results can
be deduced.  \rulex



\end{document}